\documentclass[a4paper,12pt]{article} 

\usepackage{color}

\usepackage{amsfonts, amsmath, amsthm, amssymb}
\usepackage[T1]{fontenc}
\usepackage[cp1250]{inputenc}
\usepackage{xcolor}
\usepackage{graphicx}
\usepackage{amssymb}
\usepackage{amsmath}
\usepackage{mathptmx}
\usepackage{helvet}
\usepackage{courier}
\usepackage{txfonts}
\usepackage{tikz} 
\usetikzlibrary{arrows}
\usepackage{type1cm}

\usepackage{verbatim}

\usepackage{graphicx}
\usepackage{epsfig,amscd,amssymb,amsxtra,amsmath,amsthm}
\usepackage{type1cm}
\usepackage[T1]{fontenc}
\usepackage{graphics}
\usepackage[mathscr]{eucal}
\usepackage[all]{xy}
\usepackage{amsmath,amscd}

\newcommand{\orbit}{\mathcal{O}^{\oplus}}

\newcommand{\tr}{\textup{tr}}

\newtheorem{theorem}{Theorem}[section]
\newtheorem{proposition}[theorem]{Proposition}
\newtheorem{definition}[theorem]{Definition}
\newtheorem{lemma}[theorem]{Lemma}

\newtheorem{remark}[theorem]{Remark}

\newtheorem{corollary}[theorem]{Corollary}

\newtheorem{observation}[theorem]{Observation}

\newcommand{\diam} {\mathop{\rm diam}\nolimits}

\newcommand{\Cl}  {\mathop{\rm Cl}\nolimits}

\begin{document}

\def\joinrel{\mkern-3mu}
\newcommand{\varproj}{\displaystyle \lim_{\multimapinv\joinrel-\joinrel-}}

\title{Chaos and mixing homeomorphisms on fans}
\author{Iztok Bani\v c, Goran Erceg, Judy Kennedy,  Chris Mouron and Van Nall}
\date{}

\maketitle

\begin{abstract}
We construct a mixing homeomorphism on the Lelek fan. We also construct a mixing homeomorphism on the Cantor fan. Then, we construct a family of uncountably many pairwise non-homeomorphic (non-)smooth fans  that admit a mixing homeomorphism. 
\end{abstract}
\-
\\
\noindent
{\it Keywords:} Closed relations; Mahavier products; transitive dynamical systems; transitive homeomorphisms, mixing homeomorphisms;  smooth fans, non-smooth fans\\
\noindent
{\it 2020 Mathematics Subject Classification:} 37B02, 37B45, 54C60, 54F15, 54F17

\section{Introduction}

In this paper, we study mixing homeomorphisms on compact metric spaces. By mixing, in this paper, we mean  topologically mixing. First, we study, how one can use Mahavier products of closed relations on compact metric spaces to construct a dynamical system $(X,f)$, where $f$ is a mixing homeomorphism. Then, we study quotients of dynamical systems. We start with a dynamical system $(X,f)$ and define an equivalence relation $\sim$ on $X$. Then, we discuss about when the mixing of $(X,f)$ implies  the mixing of $(X/_{\sim}, f^{\star})$. Finally, we use these techniques
\begin{enumerate}
	\item to obtain a mixing homeomorphism on the Lelek fan, 
	\item to obtain a mixing homeomorphism on the Cantor fan, and
	\item to construct a family of uncountably many pairwise non-homeomorphic (non-)smooth fans  that admit a mixing homeomorphism.
\end{enumerate} 
In addition, we show that
\begin{enumerate}
\item  there are continuous functions $f,h:L\rightarrow L$ on the Lelek fan $L$ such that
\begin{enumerate}
     \item[(a)] $h$ is a homeomorphism and $f$ is not, 
	\item[(b)]  $(L,f)$ and $(L,h)$ are both mixing as well as chaotic in the sense of Robinson but not in the sense of Devaney.
\end{enumerate}
\item  there are continuous functions $f,h:C\rightarrow C$ on the Cantor fan $C$ such that 
\begin{enumerate}
     \item $h$ is a homeomorphism and $f$ is not, 
	\item  $(C,f)$ and $(C,h)$ are both mixing as well as chaotic in the sense of  Devaney,   
\end{enumerate}
\item  there are continuous functions $f,h:C\rightarrow C$ on the Cantor fan $C$ such that 
\begin{enumerate}
     \item $h$ is a homeomorphism and $f$ is not, 
	\item  $(C,f)$ and $(C,h)$ are both mixing as well as chaotic in the sense of Robinson but not in the sense of Devaney, and
\end{enumerate}
\item  there are continuous functions $f,h:C\rightarrow C$ on the Cantor fan $C$ such that 
\begin{enumerate}
     \item $h$ is a homeomorphism and $f$ is not, 
	\item  $(C,f)$ and $(C,h)$ are both mixing as well as chaotic in the sense of Knudsen but not in the sense of Devaney.   
\end{enumerate}
\end{enumerate}

We proceed as follows. In Section \ref{s1}, we introduce the definitions, notation and the well-known results that will be used later in the paper. In Section \ref{s2}, we study mixing of Mahavier dynamical systems and mixing of  quotients of dynamical systems. Then, we use these results in Sections \ref{s3}, \ref{s4} and \ref{s5} to produce mixing homeomorphisms on various examples of fans.

\section{Definitions and Notation}\label{s1}
The following definitions, notation and well-known results are needed in the paper.


\begin{definition}
Let $X$ be a metric space, $x\in X$ and $\varepsilon>0$. We use $B(x,\varepsilon)$ to denote the open ball,  {centered} at $x$ with radius $\varepsilon$.
\end{definition}
\begin{definition}
We use $\mathbb N$ to denote the set of positive integers and $\mathbb Z$ to denote the set of integers.  
\end{definition}
\begin{definition}
Let $(X,d)$ be a compact metric space. Then we define \emph{$2^X$} by 
$$
2^{X}=\{A\subseteq X \ | \ A \textup{ is a non-empty closed subset of } X\}.
$$
Let $\varepsilon >0$ and let $A\in 2^X$. Then we define  \emph{$N_d(\varepsilon,A)$} by 
$$
N_d(\varepsilon,A)=\bigcup_{a\in A}B(a,\varepsilon).
$$
Let $A,B\in 2^X$. The function \emph{$H_d:2^X\times 2^X\rightarrow \mathbb R$}, defined by
$$
H_d(A,B)=\inf\{\varepsilon>0 \ | \ A\subseteq N_d(\varepsilon,B), B\subseteq N_d(\varepsilon,A)\},
$$
is called \emph{the Hausdorff metric}. The Hausdorff metric is in fact a metric and the metric space $(2^X,H_d)$ is called \emph{a hyperspace of the space $(X,d)$}. 
\end{definition}
\begin{remark}
Let $(X,d)$ be a compact metric space, let $A$ be a non-empty closed subset of $X$,  and let $(A_n)$ be a sequence of non-empty closed subsets of $X$. When we say $\displaystyle A=\lim_{n\to \infty}A_n$, we mean $\displaystyle A=\lim_{n\to \infty}A_n$ in $(2^X,H_d)$. 
\end{remark}
\begin{definition}
 \emph{A continuum} is a non-empty compact connected metric space.  \emph{A subcontinuum} is a subspace of a continuum, which is itself a continuum.
 \end{definition}
 \begin{definition}
Let $X$ be a continuum. 
\begin{enumerate}
\item The continuum $X$ is \emph{ unicoherent}, if for any subcontinua $A$ and $B$ of $X$ such that $X=A\cup B$,  the compactum $A\cap B$ is connected. 
\item The continuum $X$ is \emph{hereditarily unicoherent } provided that each of its subcontinua is unicoherent.
\item The continuum $X$ is a \emph{dendroid}, if it is an arcwise connected hereditarily unicoherent continuum.
\item Let $X$ be a continuum.  If $X$ is homeomorphic to $[0,1]$, then $X$ is \emph{ an arc}.   
\item A point $x$ in an arc $X$ is called \emph{an end-point of the arc  $X$}, if  there is a homeomorphism $\varphi:[0,1]\rightarrow X$ such that $\varphi(0)=x$.
\item Let $X$ be a dendroid.  A point $x\in X$ is called an \emph{end-point of the dendroid $X$}, if for  every arc $A$ in $X$ that contains $x$, $x$ is an end-point of $A$.  The set of all end-points of $X$ will be denoted by $E(X)$. 
\item A continuum $X$ is \emph{a simple triod}, if it is homeomorphic to $([-1,1]\times \{0\})\cup (\{0\}\times [0,1])$.
\item A point $x$ in a simple triod $X$ is called \emph{the top-point} or just the \emph{top of the simple triod $X$}, if  there is a homeomorphism $\varphi:([-1,1]\times \{0\})\cup (\{0\}\times [0,1])\rightarrow X$ such that $\varphi(0,0)=x$.
\item Let $X$ be a dendroid.  A point $x\in X$ is called \emph{a ramification-point of the dendroid $X$}, if there is a simple triod $T$ in $X$ with the top   $x$.  The set of all ramification-points of $X$ will be denoted by $R(X)$. 
\item The continuum $X$ is \emph{a  fan}, if it is a dendroid with at most one ramification point $v$, which is called the top of the fan $X$ (if it exists).
\item Let $X$ be a fan.   For all points $x$ and $y$ in $X$, we define  \emph{$A_X[x,y]$} to be the arc in $X$ with end-points $x$ and $y$, if $x\neq y$. If $x=y$, then we define $A_X[x,y]=\{x\}$.
\item Let $X$ be a fan with the top $v$. We say that that the fan $X$ is \emph{smooth} if for any $x\in X$ and for any sequence $(x_n)$ of points in $X$,
$$
\lim_{n\to \infty}x_n=x \Longrightarrow \lim_{n\to \infty}A_X[v,x_n]=A_X[v,x].
$$ 
\item Let $X$ be a fan.  We say that $X$ is \emph{a Cantor fan}, if $X$ is homeomorphic to the continuum
$$
\bigcup_{c\in C}{S}_c,
$$
where $C\subseteq [0,1]$ is the standard Cantor set and for each $c\in C$, ${S}_c$ is the straight line segment in the plane from $(0,0)$ to $(c,1)$. See Figure \ref{figure2}, where a Cantor fan is pictured.  
\item Let $X$ be a fan.  We say that $X$ is \emph{a Lelek fan}, if it is smooth and $\Cl(E(X))=X$. See Figure \ref{figure2}, where a Lelek fan is pictured.
\begin{figure}[h!]
	\centering
		\includegraphics[width=25em]{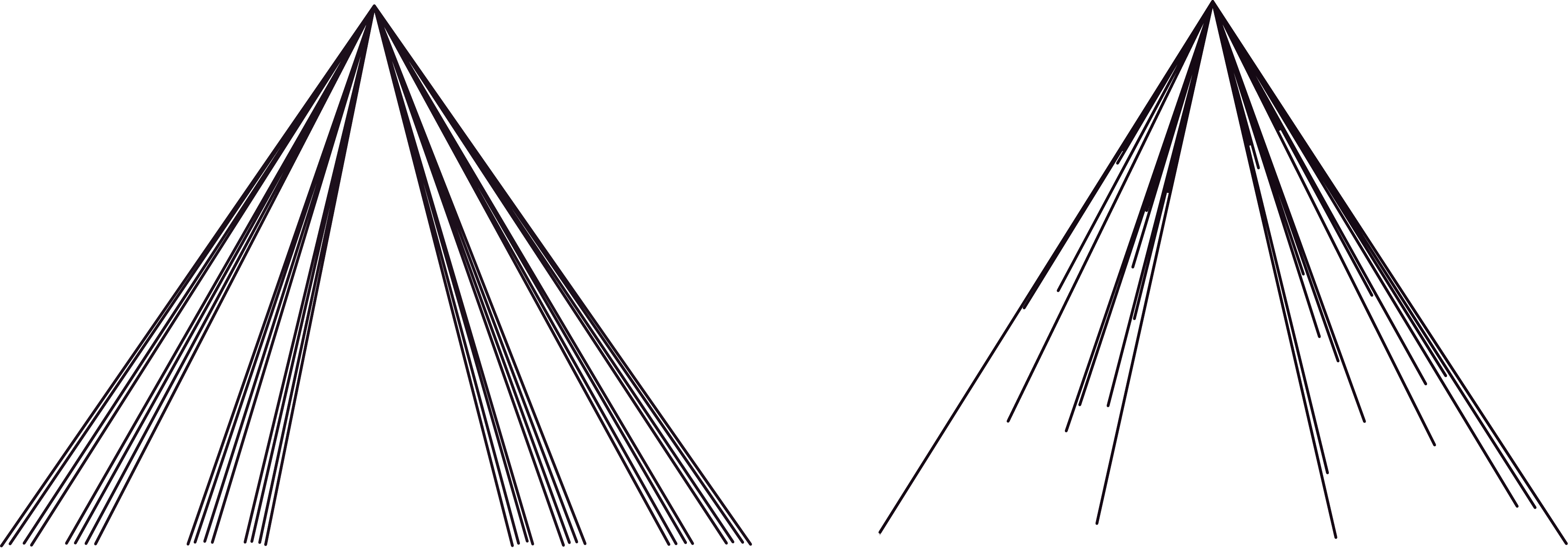}
	\caption{A Lelek fan}
	\label{figure2}
\end{figure} 
\end{enumerate}
\end{definition}
\begin{observation}\label{embi}
	It is a well-known fact that the Cantor fan is universal for smooth fans, {i.e., every smooth fan embeds into it} (for details see \cite[Theorem 9, p. 27]{Jcharatonik},  \cite[ Corollary 4]{koch},  and \cite{eberhart}).
	
	Also, note that a Lelek fan was constructed by A.~Lelek in \cite{lelek}.  An interesting  property of the Lelek fan $L$ is the fact that the set of its end-points is a dense one-dimensional set in $L$. It is also unique, i.e., any two Lelek fans are homeomorphic, for the proofs see \cite{oversteegen} and \cite{charatonik}.
\end{observation}
In this paper, $X$ will always be a non-empty compact metric space.
\begin{definition}
Let $X$ be a non-empty compact metric space and let $f:X\rightarrow X$ be a continuous function. We say that $(X,f)$ is \emph{a dynamical system}.  
\end{definition}
\begin{definition}
Let $(X,f)$ be a dynamical system and let $x\in X$. The sequence 
$$
\mathbf x=(x,f(x),f^2(x),f^3(x),\ldots)
$$   
is called \emph{the trajectory of $x$.} The set 
$$
\mathcal O_f^{\oplus}(x)=\{x,f(x),f^2(x),f^3(x),\ldots\}
$$   
is called \emph{the forward orbit set of $x$}.
\end{definition}
\begin{definition}
Let $(X,f)$ be a dynamical system and let $x\in X$.  If $\Cl(\orbit_f(x))=X$, then $x$ is called \emph{ a transitive point in $(X,f)$}.  Otherwise it is \emph{  an intransitive point in $(X,f)$}. We use \emph{  $\tr(f)$} to denote the set
$$
\tr(f)=\{x\in X \ | \ x \textup{ is a transitive point in } (X,f)\}.
$$
\end{definition}

\begin{definition}
Let $(X,f)$ be a dynamical system.  We say that $(X,f)$ is 
\emph{transitive}, if for all non-empty open sets $U$ and $V$ in $X$,  there is a non-negative integer $n$ such that $f^n(U)\cap V\neq \emptyset$. We say that the mapping $f$ is \emph{transitive}, if $(X,f)$ is transitive.
\end{definition}
The following theorem is a well-known result. See \cite{KS} for more information about transitive dynamical systems. 

\begin{theorem}\label{andrej}
Let $(X,f)$ be a dynamical system. Then the following hold.
\begin{enumerate}
	\item If $(X,f)$ is transitive, then for each $x\in \tr(f)$ and for each positive integer $n$, $f^n(x)\in \tr(f)$.
	\item If $(X,f)$ is transitive, then $\tr(f)$ is dense in $X$.
\end{enumerate}   
\end{theorem}

\begin{definition}
Let $(X,f)$ be a dynamical system.  We say that $(X,f)$ is  \emph{  mixing}, if for all non-empty open sets $U$ and $V$ in $X$,  there is a non-negative integer $n_0$ such that for each positive integer $n$,
 $$
 n\geq n_0 ~~~ \Longrightarrow ~~~  f^n(U)\cap V\neq \emptyset.
 $$
  We say that the mapping $f$ is \emph{  mixing}, if $(X,f)$ is   mixing. 
\end{definition}
\begin{definition}
	Let $(X,f)$  and $(Y,g)$ be dynamical systems. We say that 
	\begin{enumerate}
		\item $(Y,g)$ is topologically conjugate to $(X,f)$, if there is a homeomorphism $\varphi:X\rightarrow Y$ such that $\varphi\circ f=g\circ \varphi$.
		\item $(Y,g)$ is topologically semi-conjugate to $(X,f)$, if there is a continuous surjection $\alpha:X\rightarrow Y$ such that $\alpha\circ f=g\circ \alpha$. 
	\end{enumerate} 
\end{definition}

\begin{observation}\label{semi}
	Let $(X,f)$  and $(Y,g)$ be dynamical systems. Note that if $(X,f)$ is transitive (or mixing) and if $(Y,g)$ is topologically semi-conjugate to $(X,f)$, then also $(Y,g)$ is transitive (or mixing).
\end{observation}
\begin{definition}
Let $X$ be a compact metric space. We say that $X$ 
\begin{enumerate}
	\item \emph{admits a transitive homeomorphism}, if there is a homeomorphism $f:X\rightarrow X$ such that $(X,f)$ is transitive.
	\item \emph{admits a   mixing homeomorphism}, if there is a homeomorphism $f:X\rightarrow X$ such that $(X,f)$ is   mixing.
\end{enumerate} 
\end{definition}
Theorems \ref{inverz} and \ref{conju} are well-known results. Their proofs may be found in \cite{akin,aoki,KS}.
\begin{theorem}\label{inverz}
	Let $(X,f)$ be a dynamical system such that $f$ is a homeomorphism. Then the following hold. 
	\begin{enumerate}
		\item $(X,f^{-1})$ is transitive if and only if $(X,f)$ is transitive.
		\item $(X,f^{-1})$ is mixing if and only if $(X,f)$ is mixing.
	\end{enumerate} 
\end{theorem}
\begin{theorem}\label{conju}
	Let $(X,f)$ and $(Y,g)$ be dynamical systems. 
		\begin{enumerate}
		\item If $(X,f)$ is transitive and if $(Y,g)$ is topologically semi-conjugate to $(X,f)$, then $(Y,g)$ is transitive.
		\item If $(X,f)$ is mixing and if $(Y,g)$ is topologically semi-conjugate to $(X,f)$, then $(Y,g)$ is mixing.
		\end{enumerate} 
\end{theorem}
\begin{definition}
	Let $X$ be a compact metric space and let $f:X\rightarrow X$ be a continuous function. The \emph{inverse limit}, generated by $(X,f)$, is the subspace
\begin{equation*}
 \varprojlim(X,f)=\Big\{(x_{1},x_{2},x_{3},\ldots ) \in \prod_{i=1}^{\infty} X \ | \ 
\text{ for each positive integer } i,x_{i}= f(x_{i+1})\Big\}
\end{equation*}
of the topological product $\prod_{i=1}^{\infty} X$.  { The function  $\sigma : \varprojlim(X,f) \rightarrow \varprojlim(X,f)$, 
 defined by 
$$
\sigma (x_1,x_2,x_3,x_4,\ldots )=(x_2,x_3,x_4,\ldots )
$$
for each $(x_1,x_2,x_3,\ldots )\in \varprojlim(X,f)$, 
is called \emph{   the shift map on $\varprojlim(X,f)$}.    }
\end{definition}
\begin{observation}
	Note that the shift map $\sigma$ on the inverse limit $\varprojlim(X,f)$ is always a homeomorphism. Also, note that for each $(x_1,x_2,x_3,\ldots )\in \varprojlim(X,f)$,
	$$
	\sigma^{-1} (x_1,x_2,x_3,\ldots )=(f(x_1),x_1,x_2,x_3,\ldots ).
	$$
\end{observation}
Theorem Theorem \ref{shift} is a well-known result. Its proof may be found in \cite{aoki} or in \cite{KS}. 
\begin{theorem}\label{shift}
	Let $(X,f)$ be a mixing dynamical system such that $f$ is surjective and let $\sigma :\varprojlim(X,f)\rightarrow \varprojlim(X,f)$ be the shift map on $\varprojlim(X,f)$. Then the following hold.
	\begin{enumerate}
	    \item $(X,f)$ is transitive if and only if $(\varprojlim(X,f),\sigma)$ is transitive.
	    \item $(X,f)$ is mixing if and only if $(\varprojlim(X,f),\sigma)$ is mixing.
	\end{enumerate} 
\end{theorem}
\begin{definition}
		Let $(X,f)$ be a dynamical system. We say that $(X,f)$ has \emph{sensitive dependence on initial conditions}, if there is an $\varepsilon>0$ such that for each $x\in X$ and for each $\delta >0$, there are $y\in B(x,\delta)$ and a positive integer $n$ such that  
		$$
		d(f^n(x),f^n(y))>\varepsilon. 
		$$ 
			\end{definition}
			\begin{observation}\label{afna}
				Let $(X,f)$ be a dynamical system. Note that $(X,f)$ has sensitive dependence on initial conditions if and only if there is $\varepsilon >0$ such that for each non-empty open set $U$ in $X$, there is a positive integer $n$ such that $\diam (f^n(U))>\varepsilon$. See \cite[Theorem 2.22]{judyk} for more information. 
			\end{observation}
			\begin{definition}
		Let $(X,f)$ be a dynamical system and let $A$ be a non-empty closed subset of $X$. We say that \emph{$(X,f)$ has sensitive dependence on initial conditions with respect to $A$}, if there is $\varepsilon>0$ such that for each non-empty open set $U$ in $X$, there are $x,y\in U$ and a positive integer $n$ such that 
		$$
		\min\{d(f^n(x),f^n(y)),d(f^n(x),A)+d(f^n(y),A)\}>\varepsilon.
		$$  
	\end{definition}
	\begin{proposition}
		Let $(X,f)$ be a dynamical system and let $A$ be a non-empty closed subset of $X$. If $(X,f)$ has sensitive dependence on initial conditions with respect to $A$, then $(X,f)$ has sensitive dependence on initial conditions.
			\end{proposition}
			\begin{proof}
				Suppose that $(X,f)$ has sensitive dependence on initial conditions with respect to $A$ and let $\varepsilon>0$ be such that for each non-empty open set $U$ in $X$, there are $x,y\in U$ and a positive integer $n$ such that 
		$$
		\min\{d(f^n(x),f^n(y)),d(f^n(x),A)+d(f^n(y),A)\}>\varepsilon.
		$$  
        To see that  $(X,f)$ has sensitive dependence on initial conditions, we use Observation \ref{afna}. Let $U$ be any non-empty open set in $X$ and let $x,y\in U$ and let $n$ be a positive integer such that 
        $$
		\min\{d(f^n(x),f^n(y)),d(f^n(x),A)+d(f^n(y),A)\}>\varepsilon.
		$$
		Then
		$$
		\diam(f^n(U))\geq d(f^n(x),f^n(y))\geq \min\{d(f^n(x),f^n(y)),d(f^n(x),A)+d(f^n(y),A)\}>\varepsilon
		$$
		and we are done.
			\end{proof}
We use the following result.
\begin{theorem}\label{ingram}
		Let $(X,f)$ be a dynamical system, where $f$ is surjective, let $A$ be a non-empty closed subset of $X$ such that $f(A)\subseteq A$, and let
		 $\sigma$ be the shift homeomorphism on $\varprojlim(X,f)$. If $(X,f)$ has sensitive dependence on initial conditions with respect to $A$, then $(\varprojlim(X,f),\sigma^{-1})$ has sensitive dependence on initial conditions with respect to $\varprojlim(A,f|_A)$. 
	\end{theorem}
	\begin{proof}
		See \cite[Theorem 3.15]{judyk}.
	\end{proof}
We conclude this section by defining three different types of chaos. First, we define periodic points. 
\begin{definition}
		Let $(X,f)$ be a dynamical system and $p\in X$. We say that $p$ is a periodic point in $(X,f)$, if there is a positive integer $n$ such that $f^n(p)=p$. We use $\mathcal P(f)$ to denote the set of periodic points in $(X,f)$. 
	\end{definition}
	\begin{definition}
		Let $(X,f)$ be a dynamical system. We say that $(X,f)$ is chaotic in the sense of Robinson \cite{robinson}, if
		\begin{enumerate}
			\item  $(X,f)$ is transitive, and
			\item  $(X,f)$ has sensitive dependence on initial conditions.
		\end{enumerate}
	\end{definition}
	\begin{definition}
		Let $(X,f)$ be a dynamical system. We say that $(X,f)$ is chaotic in the sense of Knudsen \cite{knudsen}, if
		\begin{enumerate}
			\item  $\mathcal P(f)$ is dense in $X$, and
			\item  $(X,f)$ has sensitive dependence on initial conditions.
		\end{enumerate}
	\end{definition}
	\begin{definition}
		Let $(X,f)$ be a dynamical system. We say that $(X,f)$ is chaotic in the sense of Devaney \cite{devaney},  if
		\begin{enumerate}
			\item  $(X,f)$ is transitive, and
			\item  $\mathcal P(f)$ is dense in $X$.
		\end{enumerate}
	\end{definition}
\begin{observation}
	Note that it is proved in \cite{banks} that for any dynamical system $(X,f)$, $(X,f)$ has sensitive dependence on initial conditions, 	if $(X,f)$ is transitive and if the set $\mathcal P(f)$ is dense in $X$. 
	\end{observation}
We also use special kind of projections that are defined in the following definition. 
\begin{definition}
	For each (positive) integer $i$ and for each $\mathbf x=(x_1,x_2,x_3,\ldots)\in \prod_{k=1}^{\infty}X$ (or $\mathbf x=(\ldots,x_{-2},x_{-1},x_0,x_1,x_2,\ldots)\in \prod_{k=-\infty}^{\infty}X$ or $\mathbf x=(x_1,x_2,x_3,\ldots,x_{m})\in \prod_{k=1}^{m}X$), we use $\pi_i(\mathbf x)$ or $\mathbf x(i)$ or $\mathbf x_i$ to denote the $i$-th coordinate $x_i$ of the point $\mathbf x$.

We also use $p_1:X\times X\rightarrow X$ and $p_2:X\times X\rightarrow X$ to denote \emph{the standard projections} defined by $p_1(s,t)=s$ and $p_2(s,t)=t$ for all $(s,t)\in X\times X$.
\end{definition}

\section{Mixing, Mahavier dynamical systems and quotients of dynamical systems}
\label{s2} 
We give new results about how Mahavier products of closed relations on compact metric spaces  can be used to construct a dynamical system $(X,f)$, where $f$ is a mixing homeomorphism. Then, we study quotients of dynamical systems. Explicitly, we start with a dynamical system $(X,f)$ and an equivalence relation $\sim$ on $X$. Then, we discuss when the mixing of $(X,f)$ implies  the mixing of $(X/_{\sim}, f^{\star})$. 
\subsection{Mixing and Mahavier dynamical systems}
First, we define Mahavier products of closed relations.
\begin{definition}
Let $X$ be a {non-empty} compact metric space and let ${F}\subseteq X\times X$ be a non-empty relation on $X$. If ${F}$ is closed in $X\times X$, then we say that ${F}$ is  \emph{  a closed relation on $X$}.  
\end{definition}
\begin{definition}
Let $X$ be a {non-empty} compact metric space and let ${F}$ be a closed relation on $X$. We call
%
$$
X_F^+=\Big\{(x_1,x_2,x_3,\ldots )\in \prod_{{ i={1}}}^{\infty}X \ | \ \textup{ for each { positive} integer } i, (x_{i},x_{i+1})\in {F}\Big\}
$$
\emph{ the  Mahavier product of ${F}$}, and 
$$
X_F=\Big\{(\ldots,x_{-3},x_{-2},x_{-1},{x_0}{ ;}x_1,x_2,x_3,\ldots )\in \prod_{i={-\infty}}^{\infty}X \ | \ \textup{ for each  integer } i, (x_{i},x_{i+1})\in {F}\Big\}
$$
\emph{the two-sided  Mahavier product of ${F}$}.
\end{definition}

\begin{definition}\label{shit}
Let $X$ be a {non-empty} compact metric space and let ${F}$ be a closed relation on $X$. 
The function  $\sigma_F^{+} : {X_F^+} \rightarrow {X_F^+}$, 
 defined by 
$$
\sigma_F^{+} ({x_1,x_2,x_3,x_4},\ldots)=({x_2,x_3,x_4},\ldots)
$$
for each $({x_1,x_2,x_3,x_4},\ldots)\in {X_F^+}$, 
is called \emph{   the shift map on ${X_F^+}$}. The function  $\sigma_F : {X_F} \rightarrow {X_F}$, 
 defined by 
$$
\sigma_F (\ldots,x_{-3},x_{-2},x_{-1},{x_0};x_1,x_2,x_3,\ldots )=(\ldots,x_{-3},x_{-2},x_{-1},{x_0},x_1;x_2,x_3,\ldots )
$$
for each $(\ldots,x_{-3},x_{-2},x_{-1},{x_0};x_1,x_2,x_3,\ldots )\in {X_F}$, 
is called \emph{   the shift map on ${X_F}$}.    
\end{definition}
\begin{observation}\label{juju}
Note that $\sigma_F$ is always a homeomorphism while $\sigma_F^+$ may not be a homeomorphism.
\end{observation}
\begin{definition}
	Let $X$ be a compact metric space and let $F$ be a closed relation on $X$. The dynamical system 
	\begin{enumerate}
		\item $(X_F^{+},\sigma_F^+)$ is called \emph{a Mahavier dynamical system}.
		\item $(X_F,\sigma_F)$ is called \emph{a two-sided Mahavier dynamical system}.
	\end{enumerate}
\end{definition}
\begin{observation}\label{semiM}
	Let $X$ be a compact metric space and let $F$ be a closed relation on $X$ such that $p_1(F)=p_2(F)=X$. Note that  $(X_F^{+},\sigma_F^+)$ is semi-conjugate to $(X_F,\sigma_F)$: for $\alpha:X_F\rightarrow X_F^{+}$, $\alpha(\mathbf x)=(\mathbf x(1),\mathbf x(2), \mathbf x(3),\ldots)$ for any $\mathbf x\in X_F$, $\alpha \circ \sigma_F=\sigma_{F}^+\circ \alpha$. 
\end{observation}
 Theorems Theorem \ref{Mah} and Theorem \ref{tazadnji} are proved in \cite{BE}. We use them to prove Theorems \ref{tazadnji1} and \ref{van1}.
\begin{theorem}\label{Mah}
	Let $X$ be a compact metric space and let $F$ be a closed relation on $X$. Then
	\begin{enumerate}
		\item $\varprojlim(X_F^{+},\sigma_F^+)$ is homeomorphic to the two-sided Mahavier product $X_F$. 
		\item The inverse $\sigma_F^{-1}$ of the shift map $\sigma_F$ on $X_F$ is topologically  conjugate to the  shift map $\sigma$ on $\varprojlim(X_F^{+},\sigma_F^+)$.
	\end{enumerate}  
\end{theorem}
\begin{theorem}\label{tazadnji}
Let $X$ be a compact metric space and let $F$ be a closed relation on $X$ such that $p_1(F)=p_2(F)=X$. Then the following statements are equivalent. 
\begin{enumerate}
\item  $(X_F^+,\sigma_F^+)$ is transitive.
\item  $(X_F,\sigma_F)$ is transitive. 
\end{enumerate}
\end{theorem}
Next, we show that if  $p_1(F)=p_2(F)=X$, then  $(X_F^+,\sigma_F^+)$ is mixing if and only if $(X_F,\sigma_F)$ is mixing.
\begin{theorem}\label{tazadnji1}
Let $X$ be a compact metric space and let $F$ be a closed relation on $X$ such that $p_1(F)=p_2(F)=X$. Then the following statements are equivalent. 
\begin{enumerate}
\item  $(X_F^+,\sigma_F^+)$ is    mixing.
\item $(X_F,\sigma_F)$ is    mixing. 
\end{enumerate}
\end{theorem}
\begin{proof}
Let $\sigma$ be the shift map on $\varprojlim(X_F^{+},\sigma_F^+)$. First, suppose that  $(X_F^+,\sigma_F^+)$ is    mixing. {It follows from $p_1(F)=p_2(F)=X$ that $\sigma_F^+$ is surjective.} By  Theorem \ref{shift}, $(\varprojlim(X_F^{+},\sigma_F^+),\sigma)$ is also   mixing.  By Theorem \ref{Mah},  {$\sigma$} is   topologically conjugate to $\sigma_F^{-1}$, therefore, $(X_F,\sigma_F^{-1})$ is   mixing.  It follows from Theorem \ref{inverz} that $(X_F,\sigma_F)$ is mixing. 

Next, suppose that $(X_F,\sigma_F)$ is   mixing. By Theorem \ref{inverz}, $(X_F,\sigma_F^{-1})$ is also   mixing and it follows from Theorem \ref{Mah} that $(\varprojlim(X_F^{+},\sigma_F^+),\sigma)$ is   mixing. Since $\sigma_F^+$ is surjective, it follows from Theorem \ref{shift} that $(X_F^+,\sigma_F^+)$ is   mixing.
\end{proof}
\begin{definition}
	Let $X$ be a compact metric space. We use $\Delta_X$ to denote the diagonal-set
	$$
	\Delta_X=\{(x,x) \ | \ x\in X\}.
	$$
\end{definition}
We use the following lemma to prove Theorem \ref{van}, where we prove that for each transitive system $(X_F^+,\sigma_F^+)$, if $\Delta_X\subseteq F$, then $(X_F^+,\sigma_F^+)$ is mixing.
\begin{lemma}\label{mija1}
	Let $X$ be a compact metric space, let $F$ be a closed relation on $X$ and let $U$ be a non-empty open set in $X_F^+$. Then for each $\mathbf x\in U$, there is a positive integer $n_0$ such that for each $\mathbf y\in X_F^+$,
	$$
	\Big(\textup{for each integer } n\leq n_0,~~ \pi_n(\mathbf y)=\pi_n(\mathbf x)\Big) ~~~ \Longrightarrow ~~~ \mathbf y\in U. 
	$$
\end{lemma}
\begin{proof}
	Let $k$ be a positive integer and let $U_1$, $U_2$, $U_3$, $\ldots$, $U_k$ be open sets in $X$ such that 
	$$
	\mathbf x\in U_1\times U_2\times U_3\times \ldots \times U_k\times \prod_{i=k+1}^{\infty}X\subseteq U.
	$$
	Let $n_0=k$ and let $\mathbf y\in X_F^+$ be such that for each positive integer $n\leq n_0$, $\pi_n(\mathbf y)=\pi_n(\mathbf x)$. Then $\mathbf y\in U_1\times U_2\times U_3\times \ldots \times U_k\times \prod_{i=k+1}^{\infty}X$ and since $U_1\times U_2\times U_3\times \ldots \times U_k\times \prod_{i=k+1}^{\infty}X\subseteq U$, it follows that $\mathbf y\in U$.
\end{proof}
\begin{theorem}\label{van}
	Let $X$ be a compact metric space and let $F$ be a closed relation on $X$. If 
	\begin{enumerate}
	\item $(X_F^+,\sigma_F^+)$ is transitive, and 
	\item $\Delta_X\subseteq F$,
	\end{enumerate}
 then $(X_F^+,\sigma_F^+)$ is   mixing.
\end{theorem}
\begin{proof}
	Let $U$ and $V$ be non-empty open sets in $X_F^+$.  Since $(X_F^+,\sigma_F^+)$ is transitive, it follows from Theorem \ref{andrej} that $\tr(\sigma_F^+)$ is dense in $X_F^+$. Therefore, $\tr(\sigma_F^+)\cap U\neq \emptyset$. Let $\mathbf x\in \tr(\sigma_F^+)\cap U$. By Lemma \ref{mija1}, there is a positive integer $m_0$ such that for each $\mathbf y\in X_F^+$,
	$$
	\Big(\textup{for each positive integer } n\leq m_0,~ \pi_n(\mathbf y)=\pi_n(\mathbf x)\Big) ~~~ \Longrightarrow ~~~ \mathbf y\in U. 
	$$
 Choose and fix such a positive integer $m_0$. Next, let $m$ be a positive integer such that $m>m_0$ and such that $(\sigma_F^+)^m(\mathbf x)\in V$,  and let
 $$
 \mathbf x_1=(\mathbf x(1),\mathbf x(2),\mathbf x(3),\ldots,\mathbf x(m-1),\underbrace{\mathbf x(m),\mathbf x(m)}_{2},\mathbf x(m+1),\mathbf x(m+2),\mathbf x(m+3),\ldots).
 $$ 
 Then for each positive integer $n\leq m$, $\pi_n(\mathbf x_1)=\pi_n(\mathbf x)$. Therefore, $\mathbf x_1\in U$. Also, note that $(\sigma_F^+)^{m+1}(\mathbf x_1)=(\sigma_F^+)^{m}(\mathbf x)$, therefore, $(\sigma_F^+)^{m+1}(\mathbf x_1)\in V$.   It follows that $(\sigma_F^+)^{m+1}(U)\cap V\neq \emptyset$. 
 
 Next, let 
 $$
 \mathbf x_2=(\mathbf x(1),\mathbf x(2),\mathbf x(3),\ldots,\mathbf x(m-1),\underbrace{\mathbf x(m),\mathbf x(m),\mathbf x(m)}_{3},\mathbf x(m+1),\mathbf x(m+2),\mathbf x(m+3),\ldots).
 $$ 
 Then for each positive integer $n\leq m$, $\pi_n(\mathbf x_2)=\pi_n(\mathbf x)$. Therefore, $\mathbf x_2\in U$. Also, note that $(\sigma_F^+)^{m+2}(\mathbf x_2)=(\sigma_F^+)^{m}(\mathbf x)$, therefore, $(\sigma_F^+)^{m+2}(\mathbf x_2)\in V$.   It follows that $(\sigma_F^+)^{m+2}(U)\cap V\neq \emptyset$. 
 
 In general, let $k$ be any positive integer and let 
 $$
 \mathbf x_k=(\mathbf x(1),\mathbf x(2),\mathbf x(3),\ldots,\mathbf x(m-1),\underbrace{\mathbf x(m),\mathbf x(m),\ldots, \mathbf x(m)}_k,\mathbf x(m+1),\mathbf x(m+2),\mathbf x(m+3),\ldots).
 $$ 
 Then for each positive integer $n\leq m$, $\pi_n(\mathbf x_k)=\pi_n(\mathbf x)$. Therefore, $\mathbf x_k\in U$. Also, note that $(\sigma_F^+)^{m+k}(\mathbf x_k)=(\sigma_F^+)^{m}(\mathbf x)$, therefore, $(\sigma_F^+)^{m+k}(\mathbf x_k)\in V$.   It follows that $(\sigma_F^+)^{m+k}(U)\cap V\neq \emptyset$.  This proves that for any positive integer $n$,
 $$
 n\geq n_0 ~~~ \Longrightarrow ~~~ (\sigma_F^+)^{n}(U)\cap V\neq \emptyset,
 $$
 therefore, $(X_F^+,\sigma_F^+)$ is   mixing.
 \end{proof}
 Theorem \ref{van1} is a variant of Theorem \ref{van}, $(X_F^+,\sigma_F^+)$ from Theorem \ref{van} is replaced by $(X_F,\sigma_F)$. 
\begin{theorem}\label{van1}
	Let $X$ be a compact metric space and let $F$ be a closed relation on $X$. If 
	\begin{enumerate}
	\item $(X_F,\sigma_F)$ is transitive, and 
	\item $\Delta_X\subseteq F$,
	\end{enumerate}
 then $(X_F,\sigma_F)$ is mixing.
\end{theorem}
\begin{proof}
Suppose that $(X_F,\sigma_F)$ is transitive, and that $\Delta_X\subseteq F$. Note that $p_1(F)=p_2(F)=X$ since $\Delta_X\subseteq F$. By Theorem \ref{tazadnji}, $(X_F^+,\sigma_F^+)$ is transitive. Since $\Delta_X\subseteq F$, it follows from Theorem \ref{van} that $(X_F^+,\sigma_F^+)$ is mixing. By Theorem \ref{tazadnji1}, $(X_F,\sigma_F)$ is mixing since $\Delta_X\subseteq F$. 
\end{proof}

In Theorem \ref{judy}, we show that adding the diagonal  to the closed relation, preserves the transitivity of the Mahavier dynamical system. 
%
%

\begin{theorem}\label{judy}
	Let $X$ be a compact metric space, let $G$ be a closed relation on $X$ such that $p_1(G)=p_2(G)=X$ and let $F=G\cup \Delta_X$. Then the following hold.
	\begin{enumerate}
		\item\label{marjan} If $(X_G^+,\sigma_G^+)$ is transitive, then $(X_F^+,\sigma_F^+)$ is transitive.
		\item\label{tatjana}  If $(X_G,\sigma_G)$ is transitive, then $(X_F,\sigma_F)$ is transitive.
	\end{enumerate}
\end{theorem}
\begin{proof}
	To prove \ref{marjan}, suppose that $(X_G^+,\sigma_G^+)$ is transitive, let $m$ and $n$ be positive integers, let $U_1$, $U_2$, $U_3$, $\ldots$, $U_m$, $V_1$, $V_2$, $V_3$, $\ldots$, $V_n$ be non-empty open sets in $X$, and let 
	$$
	U=U_1\times U_2\times U_3\times \ldots\times U_m\times \prod_{k=m+1}^{\infty}X
	$$
	 and 
	 $$
	 V=V_1\times V_2\times V_3\times \ldots\times V_n\times \prod_{k=n+1}^{\infty}X
	 $$
	  be such that $U\cap X_F^+\neq \emptyset$ and $V\cap X_F^+\neq \emptyset$. To see that $(X_F^+,\sigma_F^+)$ is transitive, we prove that there is a non-negative integer $\ell$ such that $(\sigma_F^+)^{\ell}(U\cap X_F^+)\cap (V\cap X_F^+)\neq \emptyset$.

	First, let $\mathbf y\in U\cap X_F^+$ be such that $(\sigma_F^+)^{m-1}(\mathbf y)\in X_G^+$, and let 
	$$
	D=\{k \in \{1,2,3,\ldots, m-1\} \ | \ \mathbf y(k)\neq \mathbf y(k+1)\}.
	$$
	Next, let $s\in \{1,2,3,\ldots, m-1\}$ and let $k_1,k_2,k_3,\ldots,k_{s}\in \{1,2,3,\ldots ,m-1\}$ be such that
	\begin{enumerate}
		\item for each $i\in \{1,2,3,\ldots,s\}$, $k_i<k_{i+1}$ and
		\item $D=\{k_1,k_2,k_3,\ldots,k_{s}\}$.
	\end{enumerate}
	Also, let 
	\begin{align*}
	\hat U=&\underbrace{\left(\bigcap_{i=1}^{k_1}U_i\right)\times \left(\bigcap_{i=1}^{k_1}U_i\right)\times \left(\bigcap_{i=1}^{k_1}U_i\right)\times \ldots \times \left(\bigcap_{i=1}^{k_1}U_i\right)}_{k_1}\times \\
	&\underbrace{\left(\bigcap_{i=k_1+1}^{k_2}U_i\right)\times \left(\bigcap_{i=k_1+1}^{k_2}U_i\right)\times \left(\bigcap_{i=k_1+1}^{k_2}U_i\right)\times \ldots \times \left(\bigcap_{i=k_1+1}^{k_2}U_i\right)}_{k_2-k_1}\times  \ldots \\
	&\ldots \times \underbrace{\left(\bigcap_{i=k_{s-1}+1}^{k_s}U_i\right)\times \left(\bigcap_{i=k_{s-1}+1}^{k_s}U_i\right)\times \left(\bigcap_{i=k_{s-1}+1}^{k_s}U_i\right)\times \ldots \times \left(\bigcap_{i=k_{s-1}+1}^{k_s}U_i\right)}_{k_s-k_{s-1}}\times 
	\\
	&\ldots \times \underbrace{\left(\bigcap_{i=k_{s}+1}^{m}U_i\right)\times \left(\bigcap_{i=k_{s}+1}^{m}U_i\right)\times \left(\bigcap_{i=k_{s}+1}^{m}U_i\right)\times \ldots \times \left(\bigcap_{i=k_{s}+1}^{m}U_i\right)}_{m-k_{s}}\times \prod_{k=m+1}^{\infty}X
		\end{align*}
		and let 
		$$
		\bar U=\left(\bigcap_{i=1}^{k_1}U_i\right)\times \left(\bigcap_{i=k_1+1}^{k_2}U_i\right)\times \left(\bigcap_{i=k_2+1}^{k_3}U_i\right)\times \ldots \times \left(\bigcap_{i=k_{s-1}+1}^{k_s}U_i\right)\times \left(\bigcap_{i=k_s+1}^{m}U_i\right)\times \prod_{k=s+2}^{\infty}X.
		$$
Then, let $\mathbf z\in V\cap X_F^+$ be such that $(\sigma_F^+)^{n-1}(\mathbf z)\in X_G^+$, and let 
	$$
	E=\{k \in \{1,2,3,\ldots, n-1\} \ | \ \mathbf z(k)\neq \mathbf z(k+1)\}.
	$$
	Next, let $t\in \{1,2,3,\ldots, m-1\}$ and let $l_1,l_2,l_3,\ldots,l_{t}\in \{1,2,3,\ldots ,m-1\}$ be such that
	\begin{enumerate}
		\item for each $i\in \{1,2,3,\ldots,t\}$, $l_i<l_{i+1}$ and
		\item $D=\{l_1,l_2,l_3,\ldots,l_{t}\}$.
	\end{enumerate}
	Also, let 
	\begin{align*}
	\hat V=&\underbrace{\left(\bigcap_{i=1}^{l_1}V_i\right)\times \left(\bigcap_{i=1}^{l_1}V_i\right)\times \left(\bigcap_{i=1}^{l_1}V_i\right)\times \ldots \times \left(\bigcap_{i=1}^{l_1}V_i\right)}_{l_1}\times \\
	&\underbrace{\left(\bigcap_{i=l_1+1}^{l_2}V_i\right)\times \left(\bigcap_{i=l_1+1}^{l_2}V_i\right)\times \left(\bigcap_{i=l_1+1}^{l_2}V_i\right)\times \ldots \times \left(\bigcap_{i=l_1+1}^{l_2}V_i\right)}_{l_2-l_1}\times  \ldots \\
	&\ldots \times \underbrace{\left(\bigcap_{i=l_{t-1}+1}^{l_t}V_i\right)\times \left(\bigcap_{i=l_{t-1}+1}^{l_t}V_i\right)\times \left(\bigcap_{i=l_{t-1}+1}^{l_t}V_i\right)\times \ldots \times \left(\bigcap_{i=l_{t-1}+1}^{l_t}V_i\right)}_{l_t-l_{t-1}}\times 
	\\
	&\ldots \times \underbrace{\left(\bigcap_{i=l_{t}+1}^{n}V_i\right)\times \left(\bigcap_{i=l_{t}+1}^{n}V_i\right)\times \left(\bigcap_{i=l_{t}+1}^{n}V_i\right)\times \ldots \times \left(\bigcap_{i=l_{t}+1}^{n}V_i\right)}_{n-l_{t}}\times \prod_{k=n+1}^{\infty}X
		\end{align*}
		and let 
		$$
		\bar V=\left(\bigcap_{i=1}^{l_1}V_i\right)\times \left(\bigcap_{i=l_1+1}^{l_2}V_i\right)\times \left(\bigcap_{i=l_2+1}^{l_3}V_i\right)\times \ldots \times \left(\bigcap_{i=l_{t-1}+1}^{k_t}V_i\right)\times \left(\bigcap_{i=l_t+1}^{n}V_i\right)\times \prod_{k=t+2}^{\infty}X.
		$$
Note that
\begin{enumerate}
	\item $\hat U$ and $\hat V$ are both open in $\prod_{k=1}^{\infty}X$ such that 
$\mathbf y\in \hat U\subseteq U$ and $\mathbf z\in \hat V\subseteq V$, and 
    \item $\bar U$ and $\bar V$ are both open in $\prod_{k=1}^{\infty}X$ such that 
    $$
    (\mathbf y(1), \mathbf y(k_1+1), \mathbf y(k_2+1), \ldots ,\mathbf y(k_s+1), \mathbf y(k_s+2), \mathbf y(k_s+3), \ldots)\in \bar U\cap X_G^+
    $$
    and
    $$
    (\mathbf z(1), \mathbf z(l_1+1), \mathbf z(l_2+1), \ldots ,\mathbf z(l_t+1), \mathbf z(l_t+2), \mathbf z(l_t+3), \ldots)\in \bar V\cap X_G^+.
    $$
\end{enumerate} 
Next, let $\ell$ be a positive integer such that $\ell>m$ and $(\sigma_G^+)^{\ell}(\bar U\cap X_G^+)\cap (\bar V\cap X_G^+)\neq \emptyset$ and let $\bar{\mathbf x}\in \bar U\cap X_G^+$ be such that $(\sigma_G^+)^{\ell}(\bar{\mathbf x})\in  \bar V\cap X_G^+$. Note that such an integer $\ell$ does exist by Theorem \ref{andrej}. Finally, let 
\begin{align*}
	\mathbf x=&\Big(\underbrace{\bar{\mathbf x}(1),\bar{\mathbf x}(1),\bar{\mathbf x}(1),\ldots,\bar{\mathbf x}(1)}_{k_1},\underbrace{\bar{\mathbf x}(2),\bar{\mathbf x}(2),\bar{\mathbf x}(2),\ldots,\bar{\mathbf x}(2)}_{k_2-k_1},\ldots \\
	&\ldots ,\underbrace{\bar{\mathbf x}(s+1),\bar{\mathbf x}(s+1),\bar{\mathbf x}(s+1),\ldots,\bar{\mathbf x}(s+1)}_{m-k_s},\bar{\mathbf x}(s+2),\bar{\mathbf x}(s+3),\ldots, \bar{\mathbf x}(\ell), \\
	& \underbrace{\bar{\mathbf x}(\ell+1),\bar{\mathbf x}(\ell+1),\bar{\mathbf x}(\ell+1),\ldots,\bar{\mathbf x}(\ell+1)}_{l_1},\underbrace{\bar{\mathbf x}(\ell+2),\bar{\mathbf x}(\ell+2),\bar{\mathbf x}(\ell+2),\ldots,\bar{\mathbf x}(\ell+2)}_{l_2-l_1},\ldots \\
	&\ldots ,\underbrace{\bar{\mathbf x}(\ell+t+1),\bar{\mathbf x}(\ell+t+1),\bar{\mathbf x}(\ell+t+1),\ldots,\bar{\mathbf x}(\ell+t+1)}_{n-l_s},\bar{\mathbf x}(\ell+t+2),\bar{\mathbf x}(\ell+t+3),\ldots\Big)
\end{align*}
Note that $\mathbf x\in U\cap X_F^+$ and that  $\sigma_F^{\ell}(\mathbf x)\in V\cap X_F^+$. Therefore, $(\sigma_F^+)^{\ell}(U\cap X_F^+)\cap (V\cap X_F^+)\neq \emptyset$ and it follows that $(X_F^+,\sigma_F^+)$ is transitive.

To prove \ref{tatjana}, suppose that $(X_G,\sigma_G)$ is transitive. By Theorem \ref{tazadnji}, $(X_G^+,\sigma_G^+)$ is transitive, therefore, by \ref{marjan}, so is $(X_F^+,\sigma_F^+)$.  Finally, it follows from Theorem \ref{tazadnji} that $(X_F,\sigma_F)$ is transitive.
\end{proof}
\begin{corollary}\label{anamarko}
	Let $X$ be a compact metric space, let $G$ be a closed relation on $X$ such that $p_1(G)=p_2(G)=X$ and let $F=G\cup \Delta_X$. 	Then the following hold.
	\begin{enumerate}
		\item\label{uno} If $(X_G^+,\sigma_G^+)$ is transitive, then $(X_F^+,\sigma_F^+)$ is mixing.
		\item\label{due} If $(X_G,\sigma_G)$ is transitive, then $(X_F,\sigma_F)$ is mixing.
	\end{enumerate}
\end{corollary}
\begin{proof}
	To prove \ref{uno}, suppose that $(X_G^+,\sigma_G^+)$ is transitive. By Theorem \ref{judy},  $(X_F^+,\sigma_F^+)$ is transitive. Therefore, by Theorem \ref{van},  $(X_F^+,\sigma_F^+)$ is mixing since $\Delta_X\subseteq F$.
	
	To prove \ref{due}, suppose that $(X_G,\sigma_G)$ is transitive. By Theorem \ref{judy},  $(X_F,\sigma_F)$ is transitive. Therefore, by Theorem \ref{van1},  $(X_F,\sigma_F)$ is mixing since $\Delta_X\subseteq F$. 
\end{proof}

\subsection{Mixing and quotients of dynamical systems}
 Theorem \ref{tupki} is the main result of this section. First, we introduce quotients of dynamical systems and recall some of its properties.

	\begin{definition}
Let $X$ be a compact metric space and let $\sim$ be an equivalence relation on $X$. For each $x\in X$, we use $[x]$ to denote the equivalence class of the element $x$ with respect to the relation $\sim$. We also use $X/_{\sim}$ to denote the quotient space $X/_{\sim}=\{[x] \ | \ x\in X\}$. 
\end{definition}
\begin{observation}\label{kvokvo}
	Let $X$ be a compact metric space, let $\sim$ be an equivalence relation on $X$, let $q:X\rightarrow X/_{\sim}$ be the quotient map that is defined by $q(x)=[x]$ for each $x\in X$,  and let $U\subseteq X/_{\sim}$. Then 
	$$
	U \textup{ is open in } X/_{\sim} ~~~ \Longleftrightarrow ~~~ q^{-1}(U)\textup{ is open in } X.
	$$
\end{observation}
\begin{definition}\label{labod}
Let $X$ be a compact metric space, let $\sim$ be an equivalence relation on $X$,  and let $f:X\rightarrow X$ be a function such that for all $x,y\in X$,
$$
x\sim y  \Longleftrightarrow f(x)\sim f(y).
$$
 Then we let $f^{\star}:X/_{\sim}\rightarrow X/_{\sim}$ be defined by   
$
f^{\star}([x])=[f(x)]
$
for any $x\in X$. 
\end{definition}
Among other things, the following well-known proposition says that Definition \ref{labod} is good.
\begin{proposition}\label{kvocienti}
Let $X$ be a compact metric space, let $\sim$ be an equivalence relation on $X$, and  let $f:X\rightarrow X$ be a function such that for all $x,y\in X$,
$$
x\sim y  \Longleftrightarrow f(x)\sim f(y).
$$
 Then the following hold.
\begin{enumerate}
\item $f^{\star}$ is a well-defined function from  $X/_{\sim}$ to $X/_{\sim}$. 
\item If $f$ is continuous, then $f^{\star}$ is continuous.
\item If $f$ is a homeomorphism, then $f^{\star}$ is a homeomorphism.
\item If $(X,f)$ is transitive and $X/_\sim$ is metrizable, then $(X/_\sim,f^{\star})$ is transitive.
\end{enumerate}
\end{proposition}
\begin{proof}
	See \cite[Theorem 3.4]{BE}. 
\end{proof}

\begin{definition}
	Let $(X,f)$ be a dynamical system and let $\sim$ be an equivalence relation on $X$  such that for all $x,y\in X$,
$$
x\sim y  \Longleftrightarrow f(x)\sim f(y).
$$
Then we say that $(X/_{\sim},f^{\star})$ is \emph{a quotient of the dynamical system $(X,f)$} or it is \emph{the quotient of the dynamical system $(X,f)$ that is obtained from the relation $\sim$}. 
\end{definition}

\begin{observation}\label{mutula}
	Let $(X,f)$ be a dynamical system. Note that we have defined a dynamical system as a pair of a compact metric space with a continuous function on it and that in this case, $X/_{\sim}$ is not necessarily metrizable. So, if $X/_{\sim}$ is metrizable, then also $(X/_{\sim},f^{\star})$ is a dynamical system. Note that in this case, $X/_{\sim}$ is semi-conjugate to $X$: for $\alpha:X\rightarrow X/_{\sim}$, defined by $\alpha(x)=q(x)$ for any $x\in X$, where $q$ is the quotient map obtained from $\sim$, $\alpha\circ f=f^{\star}\circ \alpha$.
\end{observation}
 \begin{theorem}\label{tupki}
	Let $X$ be a compact metric space, let $\sim$ be an equivalence relation on $X$, and  let $f:X\rightarrow X$ be a function such that for all $x,y\in X$,
$$
x\sim y  \Longleftrightarrow f(x)\sim f(y).
$$
If $(X,f)$ is   mixing and $X/_\sim$ is metrizable, then $(X/_\sim,f^{\star})$ is   mixing.
\end{theorem}
\begin{proof}
	Suppose that $(X,f)$ is   mixing and that $X/_\sim$ is metrizable. It follows from Observations \ref{semi} and \ref{mutula} that  $(X/_\sim,f^{\star})$ is   mixing.
\end{proof}

\section{Mixing on the Lelek fan}\label{s3}
In this section, we produce  on the Lelek fan a mixing homeomorphism as well as a mixing mapping, which is not a homeomorphism.
\begin{definition}
	In this section, we use $X$ to denote $X=[0,1]$. For each $(r,\rho)\in (0,\infty)\times (0,\infty)$, we define the sets \emph{$L_r$}, \emph{$L_{\rho}$} and \emph{$L_{r,\rho}$}  as follows:
${L_r}=\{(x,y)\in X\times X \ | \ y=rx\}$,  ${L_{\rho}}=\{(x,y)\in X\times X \ | \ y=\rho x\}$, and $	{L_{r,\rho}}=L_r\cup L_{\rho}$.
	We also define the set \emph{$M_{r,\rho}$} by $M_{r,\rho}=X^+_{L_{r,\rho}}$.
\end{definition}
\begin{definition}
	Let {$(r,\rho)\in (0,\infty)\times (0,\infty)$}. We say that \emph{$r$ and $\rho$ never connect} or \emph{$(r,\rho)\in \mathcal{NC}$}, if \begin{enumerate}
		\item $r<1$, $\rho>1$ and 
		\item for all integers $k$ and $\ell$,  
		$$
		r^k = \rho^{\ell} \Longleftrightarrow k=\ell=0.
		$$
	\end{enumerate} 
\end{definition}

In   \cite{banic1}, the following theorem is the main result.
\begin{theorem}\label{Lelek}
Let $(r,\rho)\in \mathcal{NC}$. Then $M_{r,\rho}$ is a Lelek fan with top $(0,0,0,\ldots)$.
\end{theorem}
\begin{proof}
See \cite[Theorem 14, page 21]{banic1}.
\end{proof}
\begin{definition}\label{judydaljice}
	Let $(r,\rho)\in \mathcal{NC}$. We use $F_{r,\rho}$ to denote the following closed relation  on $X$:
	\begin{align*}
		F_{r,\rho}=L_{r,\rho}\cup\{(t,t) \ | \  t\in X\}
	\end{align*} 
	see Figure \ref{fig1}.
\begin{figure}[h!]
	\centering
		\includegraphics[width=10em]{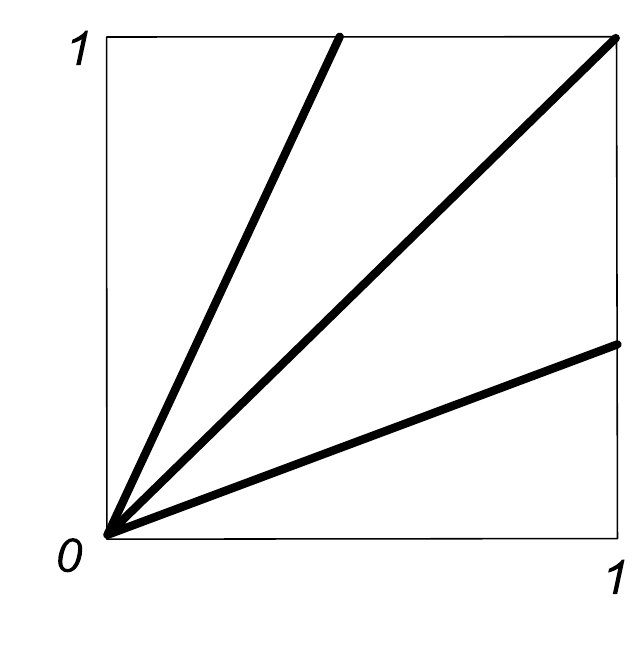}
	\caption{The relation $F$ from Definition \ref{judy}}
	\label{fig1}
\end{figure}  
\end{definition}
\begin{theorem}
	Let $(r,\rho)\in \mathcal{NC}$. Then $X_{F_{r,\rho}}^+$ and $X_{F_{r,\rho}}$ are both  Lelek fans. 
\end{theorem}
\begin{proof}
	It follows from the proof of \cite[Theorem 3.1]{EK} that $X_{F_{r,\rho}}^+$ is a Lelek fan. To see that $X_{F_{r,\rho}}$ is a Lelek fan, let 
 $$
 B_{\mathbf a,\mathbf b}=\{(\ldots,\mathbf b(2)\mathbf b(1)\cdot t,\mathbf b(1)\cdot t,t;\mathbf a(1)\cdot t, \mathbf a(2)\mathbf a(1)\cdot t,\ldots) \ | \ t\in [0,1]\}
 $$
 and
  $$
 A_{\mathbf a,\mathbf b}=B_{\mathbf a,\mathbf b}\cap X_F
 $$
 for each $\mathbf a = (\mathbf a(1), \mathbf a(2), \mathbf a(3), \ldots) \in \{1,r,\rho\} ^ {\mathbb N}$ and each $\mathbf b = (\mathbf b(1), \mathbf b(2), \mathbf b(3), \ldots) \in \{1,\frac{1}{r},\frac{1}{\rho}\} ^ {\mathbb N}$.
 Note that for each $\mathbf a\in \{1,r,\rho\}^ {\mathbb N}$ and each $\mathbf b \in \{1,\frac{1}{r},\frac{1}{\rho}\} ^ {\mathbb N}$, $B_{\mathbf a,\mathbf b}$ is a straight line segment in Hilbert cube $\prod_{k=-\infty}^{-1}[0,r^{k}]\times \prod_{k=0}^{\infty}[0,\rho^{k}]$ from $(\ldots,0,0,0;0,0,\ldots)$ to $(\ldots,\mathbf b(2)\mathbf b(1)\cdot 1,\mathbf b(1)\cdot 1,1;\mathbf a(1)\cdot 1, \mathbf a(2)\mathbf a(1)\cdot 1,\ldots)$, and that for all \emph{$\mathbf a_1,\mathbf a_2\in \{1,r,\rho\}^ {\mathbb N}$} and all $\mathbf b_1,\mathbf b_2 \in \{1,\frac{1}{r},\frac{1}{\rho}\} ^ {\mathbb N}$, 
 $$
 B_{\mathbf a_1,\mathbf b_1}\cap B_{\mathbf a_2,\mathbf b_2}=\{(\ldots,0,0,0;0,0,\ldots)\} ~~~ \Longleftrightarrow ~~~ (\mathbf a_1,\mathbf b_1)\neq (\mathbf a_2,\mathbf b_2).
 $$
Since 
$$
\Big\{(\ldots,\mathbf b(2)\mathbf b(1)\cdot 1,\mathbf b(1)\cdot 1,1;\mathbf a(1)\cdot 1, \mathbf a(2)\mathbf a(1)\cdot 1,\ldots) \  \big| \  \mathbf a\in \{1,r,\rho\}^{\mathbb N}, \mathbf b \in \big\{1,\frac{1}{r},\frac{1}{\rho}\big\} ^ {\mathbb N}\Big\}
$$
 is a Cantor set, it follows that 
 $$
 C=\bigcup_{(\mathbf{a},\mathbf{b})\in  \{1,r,\rho\}^{\mathbb N}\times \big\{1,\frac{1}{r},\frac{1}{\rho}\big\} ^ {\mathbb N}}B_{\mathbf{a},\mathbf{b}}
 $$
  is a Cantor fan.  Therefore,  $X_{F_{r,\rho}}$ is a  subcontinuum of the Cantor fan $C$.  Note that for each $\mathbf a\in \{1,r,\rho\}^{\mathbb N}$ and each $\mathbf b \in \{1,\frac{1}{r},\frac{1}{\rho}\} ^ {\mathbb N}$, $A_{\mathbf a,\mathbf b}$ is either degenerate or it is an arc from $(\ldots,0,0,0;0,0,\ldots)$ to some other point, denote it by $\mathbf e_{\mathbf a,\mathbf b}$.  Let 
$$
\mathcal U=\Big\{(\mathbf{a},\mathbf{b})\in  \{1,r,\rho\}^{\mathbb N}\times \big\{1,\frac{1}{r},\frac{1}{\rho}\big\} ^ {\mathbb N} \ | \  A_{\mathbf a,\mathbf b} \textup{ is an arc}\Big\}.
$$
Then 
$$
X_{{F_{r,\rho}}}=\bigcup_{(\mathbf{a},\mathbf{b})\in \mathcal U}A_{\mathbf a,\mathbf b}  ~~~ 
\textup{ and } ~~~ 
E(X_{F_{r,\rho}})=\{\mathbf e_{\mathbf a,\mathbf b} \ | \  (\mathbf{a},\mathbf{b})\in \mathcal U\}.
$$
Next, we show that for each $\mathbf x\in X_{F_{r,\rho}}$,
$$
\mathbf x\in E(X_{F_{r,\rho}}) ~~~ \Longleftrightarrow ~~~ \sup\{\mathbf x(k) \ | \ k \textup{ is an integer}\}=1.
$$
Let $\mathbf x\in X_{F_{r,\rho}}$.  We treat the following possible cases. 
\begin{enumerate}
\item[Case 1.] For each integer $k$, there are integers $\ell_1$ and $\ell_2$ such that $\ell_1 <k<\ell_2$ and $\mathbf x(k)\not \in \{ \mathbf x(\ell_1),\mathbf x(\ell_1)\}$.  The proof that in this case 
$$
\mathbf x\in E(X_{F_{r,\rho}}) ~~~ \Longleftrightarrow ~~~ \sup\{\mathbf x(k) \ | \ k \textup{ is an integer}\}=1
$$
follows from \cite[Theorem 3.5]{banic1}  by using the obvious homeomorphism from $X_{L_{r,\rho}}$ to the inverse limit $M=\varprojlim (M_{r,\rho},\sigma_{r,\rho})$, which is used in \cite[Section 5]{banic1} to prove that $M$ is a Lelek fan. 
 \item[Case 2.] There is an integer $k_0$ such that for each positive integer $j$, $\mathbf x (k_0-j)=\mathbf x (k_0)$ and for each integer $k$, there is an integer $\ell_0$ such that $k<\ell_0$ and $\mathbf x(k)\neq \mathbf x(\ell_0)$. The proof that  in this case
$$
\mathbf x\in E(X_{F_{r,\rho}}) ~~~ \Longleftrightarrow ~~~ \sup\{\mathbf x(k) \ | \ k \textup{ is an integer}\}=1,
$$
 is analogous to the proof of \cite[Theorem 3.5]{banic1}. 
\item[Case 3.] There is an integer $k_0$ such that for each positive integer $j$, $\mathbf x (k_0+j)=\mathbf x (k_0)$ and for each integer $k$, there is an integer $\ell_0$ such that $k>\ell_0$ and $\mathbf x(k)\neq \mathbf x(\ell_0)$. This case is analogous to the previous case.  

\item[Case 4.] There are integers $k_1$ and $k_2$ such that $k_1\leq k_2$ and such that for each positive integer $\ell$, $\mathbf x (k_1-\ell)=\mathbf x (k_1)$ and $\mathbf x (k_2+\ell)=\mathbf x (k_2)$. 
In this case, 
$$
\sup\{\mathbf x(k) \ | \ k \textup{ is an integer}\}=\max\{\mathbf x(k) \ | \ k \textup{ is an integer}\}.
$$
 Let $\mathbf x\in E(X_{F_{r,\rho}})$ and suppose that $\sup\{\mathbf x(k) \ | \ k \textup{ is an integer}\}=m<1$. Also, let $k_0$ be an integer such that $\mathbf x(k_0)=m$ and let $(\mathbf{a},\mathbf{b})\in  \{1,r,\rho\}^{\mathbb N}\times \big\{1,\frac{1}{r},\frac{1}{\rho}\big\} ^ {\mathbb N} $ be such that 
$$
\mathbf x=(\ldots,\mathbf b(2)\mathbf b(1)\cdot m,\mathbf b(1)\cdot m,m=\mathbf x(k_0),\mathbf a(1)\cdot m, \mathbf a(2)\mathbf a(1)\cdot m,\ldots ).
$$
Then 
\begin{align*}
	\mathbf x\in &\Big\{(\ldots,\mathbf b(2)\mathbf b(1)\cdot t,\mathbf b(1)\cdot t,t,\mathbf a(1)\cdot t, \mathbf a(2)\mathbf a(1)\cdot t,\ldots ) \ \big| \ t\in [0,m]\Big\}, 
\end{align*}
and 
$$
\Big\{(\ldots,\mathbf b(2)\mathbf b(1)\cdot t,\mathbf b(1)\cdot t,t,\mathbf a(1)\cdot t, \mathbf a(2)\mathbf a(1)\cdot t,\ldots ) \ \big| \ t\in [0,m]\Big\}
$$
 is a proper subarc of the arc 
\begin{align*}
	\Big\{(\ldots,\mathbf b(2)\mathbf b(1)\cdot t,\mathbf b(1)\cdot t,t,\mathbf a(1)\cdot t, \mathbf a(2)\mathbf a(1)\cdot t,\ldots ) \ \big| \ t\in [0,1]\Big\}
\end{align*}
in $X_{F_{r,\rho}}$ and is, therefore, not an endpoint of $X_{F_{r,\rho}}$. It follows that the supremum $\sup\{\mathbf x(k) \ | \ k \textup{ is an integer}\}$ equals $1$.
To prove the other implication, suppose that $\sup\{\mathbf x(k) \ | \ k \textup{ is an integer}\}=1$. Then $\mathbf x$ is the end-point of some arc $A_{\mathbf a,\mathbf b}$
in $X_{F_{r,\rho}}$, which is not equal to $(\ldots,0,0;0,0,0,\ldots)$. Therefore, it is an end-point of $X_{F_{r,\rho}}$.
\end{enumerate}
We have just proved that  
$$
\mathbf x\in E(X_{F_{r,\rho}}) ~~~ \Longleftrightarrow ~~~ \sup\{\mathbf x(k) \ | \ k \textup{ is an integer}\}=1.
$$
To see that $X_{F_{r,\rho}}$ is a Lelek fan, let $\mathbf x\in X_{F_{r,\rho}}$ be any point and let $\varepsilon >0$. We prove that there is a point $\mathbf e\in E(X_{F_{r,\rho}})$ such that $\mathbf e\in B(\mathbf x,\varepsilon)$. 
Without loss of generality, we assume that $\mathbf x\neq (\ldots, 0,0;0,0,0,\ldots)$.  Let $k_0$ be a positive integer such that $\sum_{k=k_0}^{\infty}\frac{1}{2^k}<\varepsilon$. 
It follows from \cite[Theorem 2.8]{EK} that there is a sequence $(a_1,a_2,a_3,\ldots)\in \{r,\rho\}^{\mathbf N}$ such that 
$$
\sup\{(a_1\cdot a_2\cdot a_3\cdot \ldots \cdot a_n)\cdot \mathbf x(k_0) \ | \ n \textup{ is a positive integer}\}=1.
$$
Choose and fix such a sequence $(a_1,a_2,a_3,\ldots)$.  Let
$$
\mathbf e=(\ldots,\mathbf x(-1), \mathbf x(0), \mathbf x(1), \ldots ,\mathbf x(k_0), a_1\cdot \mathbf x(k_0),a_2a_1\cdot \mathbf x(k_0),a_3a_2a_1\cdot \mathbf x(k_0),\ldots).
$$
Then $\mathbf e\in E(X_{F_{r,\rho}})$ since $\sup\{\mathbf e(k) \ | \ k \textup{ is an integer}\}=1$ and 
$$
D(\mathbf e,\mathbf x)\leq \sum_{k=k_0}^{\infty}\frac{1}{2^k}<\varepsilon,
$$
where $D$ is the metric on $X_{F_{r,\rho}}$.
This proves that also $X_{F_{r,\rho}}$ is a Lelek fan. 
\end{proof}
\begin{theorem}\label{mixLelek}
	Let $(r,\rho)\in \mathcal{NC}$. The dynamical systems $(X_{F_{r,\rho}}^+,\sigma_{F_{r,\rho}}^+)$ and $(X_{F_{r,\rho}},\sigma_{F_{r,\rho}})$ are both   mixing. 
\end{theorem}
\begin{proof}
	 It follows from \cite[Theorem 4.3 and Observation 5.3]{banic1} that $(X_{L_{r,\rho}}^+,\sigma_{L_{r,\rho}}^+)$ and  $(X_{L_{r,\rho}},\sigma_{L_{r,\rho}})$ are transitive. Since $F_{r,\rho}=L_{r,\rho}\cup \Delta_X$, it follows from Corollary \ref{anamarko} that $(X_{F_{r,\rho}}^+,\sigma_{F_{r,\rho}}^+)$ and $(X_{F_{r,\rho}},\sigma_{F_{r,\rho}})$ are both mixing.
\end{proof}
\begin{theorem}\label{robica}
The following hold for the Lelek fan $L$. 
\begin{enumerate}
\item There is a continuous mapping $f$ on the Lelek fan $L$, which is not a homeomorphism,  such that $(L,f)$ is   mixing. 
	\item There is a homeomorphism $h$ on the Lelek fan $L$ such that $(L,h)$ is   mixing. 
	\end{enumerate}
\end{theorem}
\begin{proof}
Let $(r,\rho)\in \mathcal{NC}$. We prove each part of the theorem separately. 
\begin{enumerate}
	\item Let $L=X_{F_{r,\rho}}^+$ and let $f=\sigma_F^+$. Note that $f$ is a continuous function which is not a homeomorphism. By Theorem \ref{mixLelek}, $(L,f)$ is   mixing.
	\item Let $L=X_{F_{r,\rho}}$ and let $h=\sigma_F$. Note that $h$ is a homeomorphism. By Theorem \ref{mixLelek}, $(L,h)$  is   mixing.
  \end{enumerate}
 \end{proof}
 
\section{Mixing on the Cantor fan}\label{s4}
In this section, we produce on the Cantor fan a mixing homeomorphism as well as a mixing mapping, which is not a homeomorphism. We do even more, we produce 
\begin{enumerate}
\item  continuous functions $f,h:C\rightarrow C$ on the Cantor fan $C$ such that 
\begin{enumerate}
     \item $h$ is a homeomorphism and $f$ is not, 
	\item  $(C,f)$ and $(C,h)$ are both mixing as well as chaotic in the sense of  Devaney,   
\end{enumerate}
\item  continuous functions $f,h:C\rightarrow C$ on the Cantor fan $C$ such that 
\begin{enumerate}
     \item $h$ is a homeomorphism and $f$ is not, 
	\item  $(C,f)$ and $(C,h)$ are both both mixing as well as chaotic in the sense of Robinson but not in the sense of Devaney, and
\end{enumerate}
\item  continuous functions $f,h:C\rightarrow C$ on the Cantor fan $C$ such that 
\begin{enumerate}
     \item $h$ is a homeomorphism and $f$ is not, 
	\item  $(C,f)$ and $(C,h)$ are both both mixing as well as chaotic in the sense of Knudsen but not in the sense of Devaney.   
\end{enumerate}
\end{enumerate}
We use the following theorems to prove results about periodic points.
\begin{theorem}\label{kvocki1}
	Let $(X,f)$ be a dynamical system, let $A$ be a nowhere dense closed subset of $X$ such that $f(A)\subseteq A$ and $f(X\setminus A)\subseteq X\setminus A$, and let $\sim$ be the equivalence relation on $X$, defined by 
		$$
		x\sim y ~~~  \Longleftrightarrow  ~~~ x=y \textup{ or } x,y\in A  
		$$
		for all $x,y\in X$. 
		Then the following statements are equivalent.
		\begin{enumerate}
			\item\label{ata} The set $\mathcal P(f)$ of periodic points in $(X,f)$ is dense in $X$.
			\item\label{btb} The set $\mathcal P(f^{\star})$ of periodic points in the quotient $(X/_{\sim},f^{\star})$ is dense in $X/_{\sim}$.
		\end{enumerate} 
	\end{theorem}
	\begin{proof}
		See \cite[Theorem 3.18]{judyk}.
	\end{proof}
\begin{theorem}\label{masten}
		Let $X$ be a compact metric space and let $F$ be a closed relation on $X$. If for each $(x,y)\in F$, there is a positive integer $n$ and a point $\mathbf z\in X_F^n$ such that $\mathbf z(1)=y$ and $\mathbf z(n+1)=x$, then the set of periodic points $\mathcal P(\sigma_F^+)$ is dense in $X_F^+$.
	\end{theorem}
	\begin{proof}
		See \cite[Theorem 2.18]{judyk}.
	\end{proof}
\begin{theorem}\label{ingram1}
		Let $(X,f)$ be a dynamical system and let $\sigma$ be the shift homeomorphism on $\varprojlim(X,f)$. The following statements are equivalent.
		\begin{enumerate}
			\item\label{at} The set $\mathcal P(f)$ of periodic points in $(X,f)$ is dense in $X$.
			\item\label{bt} The set $\mathcal P(\sigma^{-1})$ of periodic points in $(\varprojlim(X,f),\sigma^{-1})$ is dense in $\varprojlim(X,f)$.
		\end{enumerate} 
	\end{theorem}
	\begin{proof}
		See \cite[Theorem 3.17]{judyk}.
	\end{proof}
	We use  Theorem \ref{AS} to prove results about transitive dynamical systems on the Cantor fan.
\begin{definition}
Let $X$ be a compact metric space,  let $F$ be a closed relation on $X$ and let $x\in X$. Then we define 
$$
\mathcal U^{\oplus}_F(x)=\{y\in X \ | \ \textup{there are } n\in \mathbb N \textup{ and } \mathbf x\in X_F^{n} \textup{ such that } \mathbf x(1)=x, \mathbf x(n)=y \}
$$
and we call it the forward impression of $x$ by $F$.
\end{definition}
\begin{theorem}\label{AS}
Let $X$ be a compact metric space, let $F$ be a closed relation on $X$ and  let  $\{f_{\alpha} \ | \ \alpha \in A\}$ and $\{g_{\beta} \ | \ \beta \in B\}$ be non-empty collections of continuous functions from $X$ to $X$ such that 
$$
F^{-1}=\bigcup_{\alpha\in A}\Gamma(f_{\alpha}) ~~~  \textup{ and } ~~~ F=\bigcup_{\beta\in B}\Gamma(g_{\beta}).
$$
  If there is a dense set $D$ in $X$ such that for each $s\in D$,  $\Cl(\mathcal U^{\oplus}_F(s))=X$, then $(X_F^+,\sigma_F^+)$ is transitive. 
\end{theorem}
\begin{proof}
	See \cite[Theorem 4.8]{BE}.
\end{proof}Finally, we use the following theorem when studying sensitive dependence on initial conditions.
\begin{theorem}\label{kvocki}
	Let $(X,f)$ be a dynamical system, let $A$ be a nowhere dense closed subset of $X$ such that $f(A)\subseteq A$ and $f(X\setminus A)\subseteq X\setminus A$, and let $\sim$ be the equivalence relation on $X$, defined by 
		$$
		x\sim y ~~~  \Longleftrightarrow  ~~~ x=y \textup{ or } x,y\in A  
		$$
		for all $x,y\in X$. 
The following statements are equivalent.
	\begin{enumerate}
		\item\label{uno1} $(X,f)$ has sensitive dependence on initial conditions with respect to $A$. 
		\item\label{due2} $(X/_{\sim},f^{\star})$ has sensitive dependence on initial conditions.
	\end{enumerate}	 
	\end{theorem}
	\begin{proof}
		See \cite[Theorem 3.16]{judyk}.
	\end{proof}
\subsection{Mixing and Devaney's chaos on the Cantor fan}\label{31}
  Here, we study functions $f$ on the Cantor fan $C$ such that $(C,f)$ is mixing as we well as chaotic in the sense of Devaney.

  \begin{definition}\label{nall2}
In this subsection, we use $X$ to denote $X=[0,1]\cup [2,3]\cup [4,5]\cup [6,7]\cup [8,9]$, and we let $f_1,f_2,f_3:X\rightarrow X$ to be the homeomorphisms from $X$ to $X$ that are defined by  
$$
f_1(x)=x,
$$
$$
f_2(x)=\begin{cases}
                x\text{;} & x\in [8,9]\\
				x+2\text{;} & x\in [0,1]\cup[4,5]\\
				(x-2)^{2}\text{;} & x\in [2,3]\\
				(x-6)^{3}+4\text{;} & x\in [6,7]
			\end{cases} 
$$
$$
f_3(x)=\begin{cases}
				x\text{;} & x\in [0,1]\\
				x+2\text{;} & x\in [2,3]\cup[6,7]\\
				(x-4)^{\frac{1}{2}}+2\text{;} & x\in [4,5]\\
				(x-8)^{\frac{1}{3}}+6\text{;} & x\in [8,9]
			\end{cases}  
$$
for each $x\in X$. Then we use $F$ to denote the relation 
$$
F=\Gamma(f_1)\cup \Gamma(f_2)\cup \Gamma(f_3);
$$
 see Figure \ref{figca2}.
\begin{figure}[h!]
	\centering
		\includegraphics[width=30em]{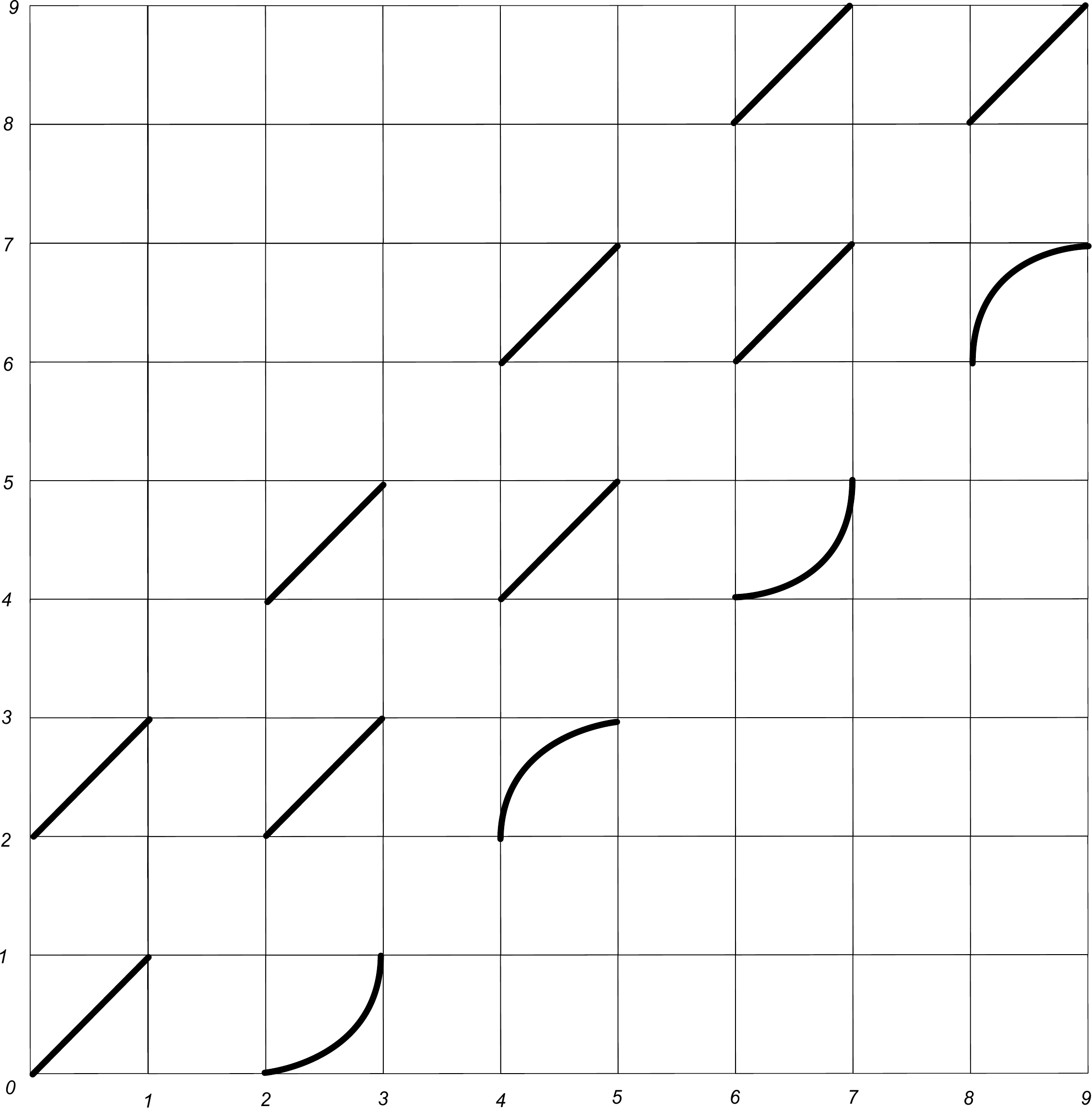}
	\caption{The relation $F$ from Definition \ref{nall2}}
	\label{figca2}
\end{figure}  
	\end{definition}
\begin{definition}
	We define two equivalence relations. 
	\begin{enumerate}
		\item For all $\mathbf x,\mathbf y\in X_F^+$, we define the relation $\sim_+$ as follows:
$$
\mathbf x\sim_+ \mathbf y ~~~   \Longleftrightarrow  ~~~  \mathbf x = \mathbf y \textup{ or  for each positive integer } k,  {\{\mathbf x(k),\mathbf y(k)\}\subseteq} \{0,2,4,6,8\}. 
$$
\item For all $\mathbf x,\mathbf y\in X_F$, we define the relation $\sim$ as follows:
$$
\mathbf x\sim \mathbf y ~~~   \Longleftrightarrow  ~~~  \mathbf x = \mathbf y \textup{ or  for each integer } k,  {\{\mathbf x(k),\mathbf y(k)\}\subseteq} \{0,2,4,6,8\}. 
$$
	\end{enumerate}
\end{definition}
\begin{observation}
	Essentially the same proof as the one from \cite[Example 4.14]{BE} shows that the quotient spaces $X_F^+{/}_{\sim_+}$ and $X_F{/}_{\sim}$ are both  Cantor fans. Also, note that $(\sigma_F^+)^{\star}$ is not a homeomorphism on $X_F^+{/}_{\sim_+}$ while $\sigma_F^{\star}$ is a homeomorphism on $X_F{/}_{\sim}$.
	\end{observation}
	
		\begin{theorem}\label{ppp}
The following hold for the sets of periodic points in $(X_F^+/_{\sim_+},(\sigma_F^+)^{\star})$ and $(X_F/_{\sim},\sigma_F^{\star})$.
\begin{enumerate}
	\item\label{lubi} The set $\mathcal P((\sigma_F^+)^{\star})$ of periodic points in the quotient $(X_F^+/_{\sim_+},(\sigma_F^+)^{\star})$ is dense in $X_F^+/_{\sim_+}$. 
	\item The set $\mathcal P(\sigma_F^{\star})$ of periodic points in the quotient $(X_F/_{\sim},\sigma_F^{\star})$ is dense in $X_F/_{\sim}$. 
\end{enumerate}
	\end{theorem}
\begin{proof}
Using Theorem \ref{kvocki1}, we prove each of the statements separately.
\begin{enumerate}
	\item We use Theorem \ref{masten} to prove the first part of the theorem. Let $(x,y)\in F$ be any point. We show that there are a positive integer $n$ and a point $\mathbf z\in X_F^n$ such that $\mathbf z(1)=y$ and $\mathbf z(n+1)=x$. We consider the following cases for $x$.
	\begin{enumerate}
	\item $x\in [0,1]$. If $y=x$, then let $n=1$ and $\mathbf z=(x,x)$. If $y=x+2$, then let $n=1$ and $\mathbf z=(x+2,x)$. 
	\item $x\in [2,3]$. If $y=(x-2)^2$, then let $n=3$ and $\mathbf z=((x-2)^2,(x-2)^2+2,(x-2)^2+4,x)$. If $y=x$, then let $n=1$ and $\mathbf z=(x,x)$. If $y=x+2$, then let $n=3$ and $\mathbf z=(x+2,(x-2)^{\frac{1}{2}}+2,x-2,x)$. 
	\item $x\in [4,5]$. If $y=(x-4)^{\frac{1}{2}}+2$, then let $n=3$ and $\mathbf z=((x-4)^{\frac{1}{2}}+2,x-4,x-2,x)$.  If $y=x$, then let $n=1$ and $\mathbf z=(x,x)$. If $y=x+2$, then let $n=3$ and $\mathbf z=(x+2,x+4,(x-4)^{\frac{1}{3}}+6, x)$. 
	\item $x\in [6,7]$. If $y=(x-6)^{3}+4$, then let $n=3$ and $\mathbf z=((x-6)^{3}+4,(x-6)^{3}+6,(x-6)^{3}+8,x)$. If $y=x$, then let $n=1$ and $\mathbf z=(x,x)$. If $y=x+2$, then let $n=3$ and $\mathbf z=(x+2,(x-6)^{\frac{1}{3}}+6,x-2,x)$.
	\item $x\in [8,9]$. If $y=(x-8)^{\frac{1}{3}}+6$, then let $n=3$ and $\mathbf z=((x-8)^{\frac{1}{3}}+6,x-4,x-2,x)$. If $y=x$, then let $n=1$ and $\mathbf z=(x,x)$. 
	\end{enumerate}	
	\item It follows from \ref{lubi} and from Theorem \ref{ingram1} that the set $\mathcal P(\sigma^{-1})$ of periodic points in $(\varprojlim(X_F^+,\sigma_F^+),\sigma^{-1})$ is dense in $\varprojlim(X_F^+,\sigma_F^+)$. By Theorem \ref{Mah}, the set $\mathcal P(\sigma_F^{\star})$ of periodic points in the quotient $(X_F/_{\sim},\sigma_F^{\star})$ is dense in $X_F/_{\sim}$. 
\end{enumerate}
\end{proof}

\begin{theorem}\label{ttt}
	The dynamical systems $(X_F^+/_{\sim_+},(\sigma_F^+)^{\star})$ and $(X_F/_{\sim},\sigma_F^{\star})$ are both transitive.
	\end{theorem}
\begin{proof}
	To prove that $(X_F^+/_{\sim_+},(\sigma_F^+)^{\star})$ is transitive, we prove that $(X_F^+,\sigma_F^+)$ is transitive.  Note that both $F$ and $F^{-1}$ are unions of three graphs of homeomorphisms. So, all the initial conditions from Theorem \ref{AS} are satisfied. To see that $(X_F^+,\sigma_F^+)$ is transitive, we prove that there is a dense set $D$ in $X$ such that for each $s\in D$,  $\Cl(\mathcal U^{\oplus}_H(s))=X$. Let $D=(0,1)\cup(2,3)\cup(4,5)\cup(6,7)\cup(8,9)$. Then $D$ is dense in $X$. Let $s\in D$ be any point and let $\ell\in\{0,1,2,3,4\}$ be such that $s\in (2\ell,2\ell+1)$. Note that 
	$$
	s,s-2,s-4, s-6,\ldots, s-2\ell\in \mathcal U^{\oplus}_H(s)
	$$
	and let $t=s-2\ell$. Then $t\in (0,1)$. It follows from the definition of $F$ that for all integers $m$, $n$ and for each $k\in\{0,1,2,3,4\}$, 
	$$
	t^{\frac{2^m}{3^n}}+k\cdot 2\in \mathcal U^{\oplus}_F(t).
	$$
	It follows from Theorem \cite[Lemma 4.9]{BE} that $\big\{t^{\frac{2^m}{3^n}}+k\cdot 2 \ | \ m,n\in \mathbb Z, k\in \{0,1,2,3,4\}\big\}$ is dense in $X$. Since 
	$$
	\big\{t^{\frac{2^m}{3^n}}+k\cdot 2 \ | \ m,n\in \mathbb Z, k\in \{0,1,2,3,4\}\big\}\subseteq \mathcal U^{\oplus}_F(t)\subseteq \mathcal U^{\oplus}_F(s),
	$$
	it follows that $\mathcal U^{\oplus}_F(s)$ is dense in $X$. Therefore, by Theorem \ref{AS}, $(X_F^+,\sigma_F^+)$ is transitive and it follows from Theorem \ref{tazadnji} that $(X_F,\sigma_F)$ is transitive since $p_1(F)=p_2(F)=X$. It follows from Theorem \ref{kvocienti} that $(X_F^+/_{\sim_+},(\sigma_F^+)^{\star})$ and $(X_F/_{\sim},\sigma_F^{\star})$ are both transitive. 
\end{proof}
\begin{theorem}\label{sss}
		The dynamical systems $(X_F^+/_{\sim_+},(\sigma_F^+)^{\star})$ and $(X_F/_{\sim},\sigma_F^{\star})$ both have sensitive dependence on initial conditions. 
	\end{theorem}
\begin{proof}
The dynamical systems $(X_F^+/_{\sim_+},(\sigma_F^+)^{\star})$ and $(X_F/_{\sim},\sigma_F^{\star})$ are both transitive by Theorem \ref{ttt}. Also, by Theorem \ref{ppp}, the set $\mathcal P((\sigma_F^+)^{\star})$ of periodic points in the quotient $(X_F^+/_{\sim_+},(\sigma_F^+)^{\star})$ is dense in $X_F^+/_{\sim_+}$,  and the set $\mathcal P(\sigma_F^{\star})$ of periodic points in the quotient $(X_F/_{\sim},\sigma_F^{\star})$ is dense in $X_F/_{\sim}$. It follows from  \cite[Theorem]{banks} that $(X_F^+/_{\sim_+},(\sigma_F^+)^{\star})$ and $(X_F/_{\sim},\sigma_F^{\star})$ both have sensitive dependence on initial conditions.
\end{proof}

\begin{theorem}\label{robi}
The following hold for the Cantor fan $C$. 
\begin{enumerate}
\item There is a continuous mapping $f$ on the Cantor fan $C$, which is not a homeomorphism,  such that $(C,f)$ is  mixing as well as chaotic in the sense of  Devaney.
	\item There is a homeomorphism $h$ on the Cantor fan $C$ such that $(C,h)$ is  mixing as well as chaotic in the sense of Devaney. 
	\end{enumerate}
\end{theorem}
\begin{proof}
We prove each part of the theorem separately. 
\begin{enumerate}
	\item Let $C=X_F^+/_{\sim_+}$ and let $f=(\sigma_F^+)^{\star}$. Note that $f$ is a continuous function which is not a homeomorphism. By Theorem \ref{sss}, $(C,f)$ has sensitive dependence on initial conditions, by Theorem \ref{ttt}, $(C,f)$ is transitive, and by Theorem \ref{ppp}, the set $\mathcal P(f)$ of periodic points in $(C,f)$ is dense in $C$.  Therefore, $(C,f)$ is chaotic in the sense of Devaney. 
	
	It follows from Theorem \ref{van} that $(X_F^+,\sigma_F^+)$ is mixing since $\Delta_X\subseteq F$. It follows from Theorem \ref{tupki} that $(C,f)$ is also mixing.
	\item Let $C=X_F/_{\sim}$ and let $h=\sigma_F^{\star}$. Note that $h$ is a homeomorphism. By Theorem \ref{sss}, $(C,h)$ has sensitive dependence on initial conditions, by Theorem \ref{ttt}, $(C,h)$ is transitive, and by Theorem \ref{ppp}, the set $\mathcal P(h)$ of periodic points in  $(C,h)$ is dense in $C$.  Therefore, $(C,h)$ is chaotic in the sense of Devaney. 
	
	It follows from Theorem \ref{van1} that $(X_F,\sigma_F)$ is mixing since $\Delta_X\subseteq F$. It follows from Theorem \ref{tupki} that $(C,h)$ is also mixing.
  \end{enumerate}
 \end{proof}

 \subsection{Mixing and Robinson's but not Devaney's chaos on the Cantor fan}\label{32} 
  Here, we study functions $f$ on the Cantor fan $C$ such that $(C,f)$ is mixing as we well as chaotic in the sense of Robinson but not in the sense of Devaney.

	\begin{definition}\label{nall1}
In this subsection, we use $X$ to denote 
$$
X=[0,1]\cup [2,3]\cup[4,5]
$$
 and we let $f_1,f_2,f_3:X\rightarrow X$ to be the homeomorphisms from $X$ to $X$ that are defined by  
$$
f_1(x)=x,
$$
$$
f_2(x)=\begin{cases}
				x+2\text{;} & x\in [0,1]\\
				(x-2)^2\text{;} & x\in [2,3]\\
				x\text{;} & x\in [4,5]
			\end{cases} 
$$
$$
f_3(x)=\begin{cases}
				x\text{;} & x\in [0,1]\\
				x+2\text{;} & x\in [2,3]\\
				(x-4)^{\frac{1}{3}}+2\text{;} & x\in [4,5]			
				\end{cases}  
$$
for each $x\in X$. Then we use $F$ to denote the relation 
$$
F=\Gamma(f_1)\cup \Gamma(f_2)\cup \Gamma(f_3);
$$
 see Figure \ref{fig2}.
\begin{figure}[h!]
	\centering
		\includegraphics[width=14em]{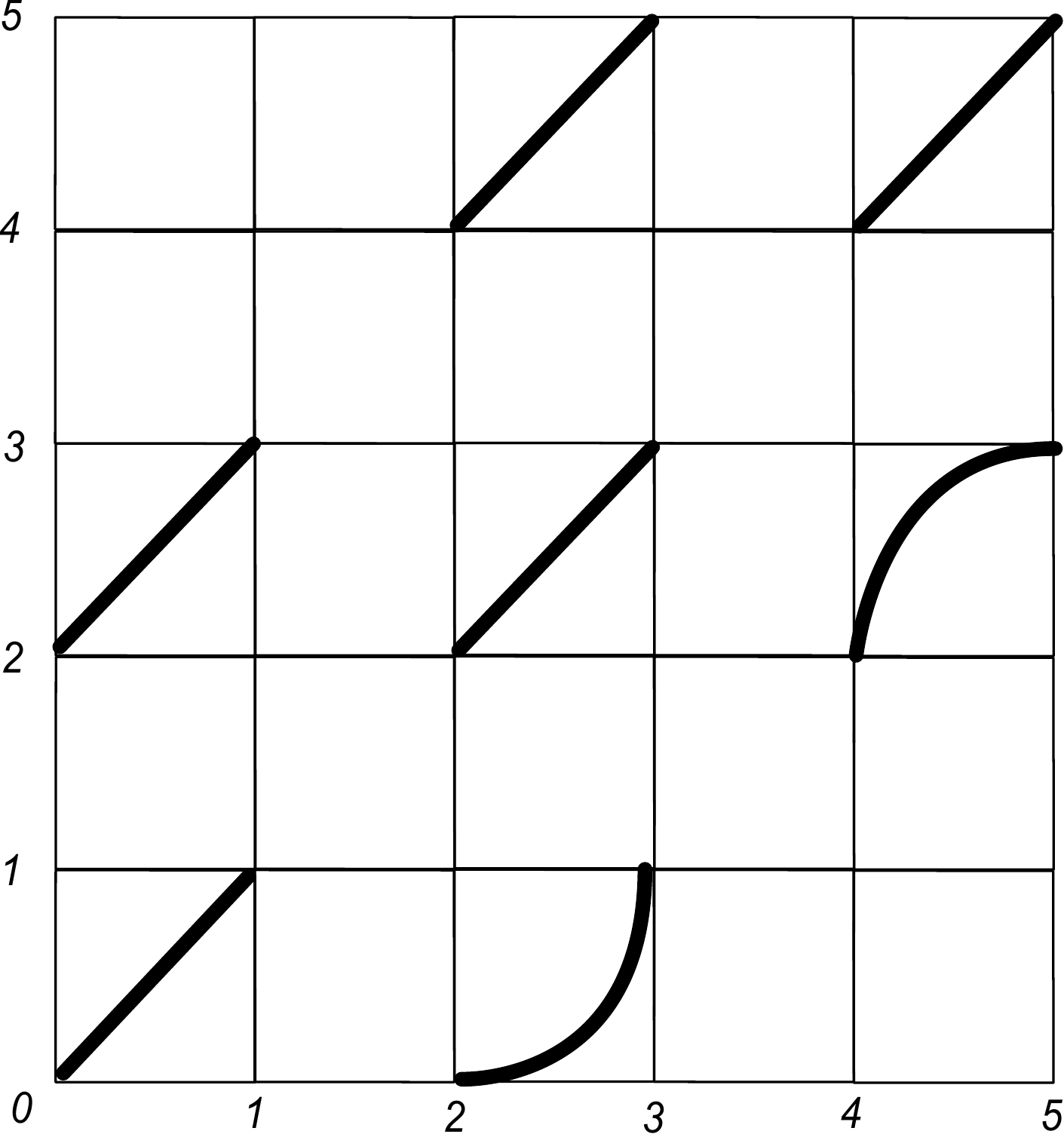}
	\caption{The  relation $F$ from Definition \ref{nall1}}
	\label{fig2}
\end{figure}  
	\end{definition}
	\begin{definition}
	We define two equivalence relations. 
	\begin{enumerate}
		\item For all $\mathbf x,\mathbf y\in X_F^+$, we define the relation $\sim_+$ as follows:
$$
\mathbf x\sim_+ \mathbf y ~~~   \Longleftrightarrow  ~~~  \mathbf x = \mathbf y \textup{ or  for each positive integer } k,  {\{\mathbf x(k),\mathbf y(k)\}\subseteq} \{0,2,4\}. 
$$
\item For all $\mathbf x,\mathbf y\in X_F$, we define the relation $\sim$ as follows:
$$
\mathbf x\sim \mathbf y ~~~   \Longleftrightarrow  ~~~  \mathbf x = \mathbf y \textup{ or  for each integer } k,  {\{\mathbf x(k),\mathbf y(k)\}\subseteq} \{0,2,4\}. 
$$
	\end{enumerate}
\end{definition}
\begin{observation}
	Note that it follows from \cite[Example 4.14]{BE} that the quotient spaces $X_F^+{/}_{\sim_+}$ and $X_F{/}_{\sim}$ are both  Cantor fans. Also, note that $(\sigma_F^+)^{\star}$ is not a homeomorphism on $X_F^+{/}_{\sim_+}$ while $\sigma_F^{\star}$ is a homeomorphism on $X_F{/}_{\sim}$.5
	\end{observation}
	First, we prove the following theorems about sensitive dependence on initial conditions.
	\begin{theorem}\label{skola0}
		Let $A=\big\{\mathbf x\in X_F^+ \ | \ \textup{ for each positive integer } k, ~~ \mathbf x(k)\in \{0,2,4\}\big\}$.
		 Then
		 \begin{enumerate}
		 	\item $\sigma_F^+(A)\subseteq A$ and $\sigma_F^+(X_F^+\setminus A)\subseteq X_F^+\setminus A$, and
		 	\item $(X_F^+,\sigma_F^+)$ has sensitive dependence on initial conditions with respect to $A$.
		 \end{enumerate}  
	\end{theorem}
	\begin{proof}
		First, note that $\sigma_F^+(A)\subseteq A$ and $\sigma_F^+(X_F^+\setminus A)\subseteq X_F^+\setminus A$. Next, let $f=\sigma_F^+$ and let $\varepsilon =\frac{1}{4}$. We show that for each basic open set $U$ of the product topology on $\prod_{k=1}^{\infty}X$ such that $U\cap X_F^+\neq \emptyset$, there are $\mathbf x,\mathbf y\in U\cap X_F^+$ such that for some positive integer $m$,
$$
\min\{d(f^m(\mathbf x),f^m(\mathbf y)),d(f^m(\mathbf x),A)+d(f^m(\mathbf y),A)\}>\varepsilon,
$$
 where $d$ is the product metric on $\prod_{k=1}^{\infty}X$, defined by  	
		$$
		d((x_1,x_2,x_3,\ldots),(y_1,y_2,y_3,\ldots))=\max\Big\{\frac{|y_k-x_k|}{2^k} \ \big| \  k \textup{ is a positive integer}\Big\}
		$$ 
		for all $(x_1,x_2,x_3,\ldots),(x_1,x_2,x_3,\ldots)\in \prod_{k=1}^{\infty}X$.
  Let $U$ be a basic set of the product topology on $\prod_{k=1}^{\infty}X$ such that $U\cap X_F^+\neq \emptyset$. Also, let $n$ be a positive integer and for each $i\in \{1,2,3,\ldots, n\}$, let $U_i$ be an open set in $X$ such that 
	$$
	U=U_1\times U_2\times U_3\times \ldots \times U_n\times \prod_{k=n+1}^{\infty}X. 
	$$
Next, let $\mathbf z=(z_1,z_2,z_3,\ldots)\in U\cap X_F^+$ be any point such that $z_n\not \in \{0,1,2,3,4,5\}$. 
We consider the following possible cases for the coordinate $z_n$ of the point $\mathbf z$.
\begin{enumerate}
\item $z_n\in (0,1)$. Then let $\mathbf x=(x_1,x_2,x_3,\ldots)\in X_F^+$ be defined by 
	$$
	(x_1,x_2,x_3,\ldots,x_n)=(z_1,z_2,z_3,\ldots,z_n)
	$$
and for each positive integer $k$, $x_{n+k}=z_n$. Also, we define $\mathbf y=(y_1,y_2,y_3,\ldots)\in X_F^+$ as follows. First, let 
$$
	(y_1,y_2,y_3,\ldots,y_n)=(z_1,z_2,z_3,\ldots,z_n).
	$$
Next, we define 
\begin{align*}
	&(y_{n+1},y_{n+2},y_{n+3},\ldots)=\\
	&\Big(z_n+2,z_n+4,z_n^{\frac{1}{3}}+2, z_n^{\frac{1}{3}}+4,z_n^{\frac{1}{3^2}}+2,z_n^{\frac{1}{3^2}}+4,z_n^{\frac{1}{3^3}}+2,z_n^{\frac{1}{3^3}}+4,\ldots\Big).
\end{align*}
	Note that 
	$$
	\lim_{k\to \infty}y_{n+4+2k}=5 ~~~ \textup{ and } ~~~ \lim_{k\to \infty}y_{n+3+2k}=3.
	$$
	Let $k_0$ be an even positive integer such that for each positive integer $k$,
	$$
	k\geq k_0 ~~~ \Longrightarrow ~~~  5-y_{n+4+2k}<\frac{1}{10} \textup{ and } 3-y_{n+3+2k}<\frac{1}{10}.
	$$

	Let $m=n+k_0+1$. Then,
	\begin{align*}
		d(f^m(\mathbf x),f^m(\mathbf y))=&\max\Big\{\frac{|y_k-x_k|}{2^{k-m+1}} \ \big| \  k \in \{m,m+1,m+2,m+3,\ldots\}\Big\}\geq 1>\varepsilon
	\end{align*}
	and 
	\begin{align*}
		d(f^m(\mathbf x),A)+d(f^m(\mathbf y),A)\geq &d(f^m(\mathbf y),A)=\min\{d(f^m(\mathbf y),\mathbf a) \ | \ \mathbf a\in A\}	=\\
		&\min \Bigg\{\max\Big\{\frac{|\mathbf a(k)-y_{k+m}|}{2^{k}} \ \big| \  k \in \{1,2,3,\ldots\}\Big\} \ \Big| \ \mathbf a\in A\Bigg\}\geq \\
		& 	\frac{\mathbf y_{k+m}-4}{2}\geq \frac{9}{20}	>\varepsilon \\
			\end{align*}
	\item $z_n \not\in (0,1)$. Then there is an integer $j\in \{1,2,3\}$ such that $z_n\in (2j,2j+1)$. In this case, the proof is analogous to the proof of the previous case. We leave the details to the reader.
	    \end{enumerate}
    This proves that $(X_F^+,\sigma_F^+)$ has sensitive dependence on initial conditions with respect to $A$.
	\end{proof}	
	\begin{corollary}\label{corolar}
	Let $B=\big\{\mathbf x\in X_F \ | \ \textup{ for each integer } k, \mathbf x(k)\in \{0,2,4\}\big\}$.
		 Then
		 \begin{enumerate}
		 	\item $\sigma_F(B)\subseteq B$ and $\sigma_F(X_F\setminus B)\subseteq X_F\setminus B$, and
		 	\item $(X_F,\sigma_F)$ has sensitive dependence on initial conditions with respect to $B$.
		 \end{enumerate}
	\end{corollary}
	\begin{proof}
	First, note that $\sigma_F(B)\subseteq B$ and $\sigma_F(X_F\setminus B)\subseteq X_F\setminus B$. Next, let 
	$$
	A=\big\{\mathbf x\in X_F^+ \ | \ \textup{ for each positive integer } k, \mathbf x(k)\in \{0,2,4,6\}\big\}.
	$$
	 By Theorem \ref{skola0},
		 \begin{enumerate}
		 	\item $\sigma_F^+(A)\subseteq A$ and $\sigma_F^+(X_F^+\setminus A)\subseteq X_F^+\setminus A$, and
		 	\item $(X_F^+,\sigma_F^+)$ has sensitive dependence on initial conditions with respect to $A$.
		 \end{enumerate}
		 Note that $\sigma_F^+$ is surjective. By Theorem \ref{ingram}, $		  (\varprojlim(X_F^+,\sigma_F^+),\sigma^{-1})$
		  has sensitive dependence on initial conditions with respect to $\varprojlim(A,\sigma_F^+|_{A})$, where $\sigma$ is the shift homeomorphism on $\varprojlim(X_F^+,\sigma_F^+)$. By Theorem \ref{Mah}, the inverse limit $\varprojlim(X_F^+,\sigma_F^+)$ is homeomorphic to the two-sided Mahavier product $X_F$ and the inverse of the shift homeomorphism $\sigma_F$ on $X_F$ is topologically conjugate to the  shift homeomorphism $\sigma $ on $\varprojlim(X_F^{+},\sigma_F^+)$. Let $ \varphi:\varprojlim(X_F^{+},\sigma_F^+)\rightarrow X_F$ be the homeomorphism, used to prove Theorem \ref{Mah} in \cite[Theorem 4.1]{BE}. Then $\varphi(\varprojlim(A,\sigma_F^+|_{A}))=B$.
 		   	 Therefore, $(X_F,\sigma_F)$ has sensitive dependence on initial conditions with respect to $B$. 
 		   	  \end{proof}
	
	\begin{theorem}\label{cantor}
		The dynamical systems $(X_F^+/_{\sim_+},(\sigma_F^+)^{\star})$ and $(X_F/_{\sim},\sigma_F^{\star})$ both have sensitive dependence on initial conditions. 
	\end{theorem}
\begin{proof}
For each of the dynamical systems, we prove separately  that it has sensitive dependence on initial conditions.
\begin{enumerate}
	\item Let $C=X_F^+/_{\sim_+}$ and let $f=(\sigma_F^+)^{\star}$, i.e., for each $\mathbf x\in X_F$, $f([\mathbf x])=[\sigma_F^+(\mathbf x)]$.
We show that $(C,f)$ has sensitive dependence on initial conditions. 
Let 
$$
A=\big\{\mathbf x\in X_F^+ \ | \ \textup{ for each positive integer } k, \mathbf x(k)\in \{0,2,4\}\big\}.
$$
		 By Theorem  \ref{skola0}, 
		 \begin{enumerate}
		 	\item $\sigma_F(A)\subseteq A$ and $\sigma_F(X_F^+\setminus A)\subseteq X_F^+\setminus A$ and
		 	\item $(X_F^+,\sigma_F^+)$ has sensitive dependence on initial conditions with respect to $A$.
		 \end{enumerate}
Since $A$ is a closed nowhere dense set in $X_F^+$, it follows from Theorem \ref{kvocki} that $(C,f)$ has sensitive dependence on initial conditions.
\item Let $C=X_F/_{\sim}$ and let $h=\sigma_F^{\star}$, i.e., for each $\mathbf x\in X_F$, $h([\mathbf x])=[\sigma_F(\mathbf x)]$.
We show that $(C,h)$ has sensitive dependence on initial conditions. The rest of the proof is analogous to the proof above - instead of the set $A$, the set  
$$
B=\big\{\mathbf x\in X_F \ | \ \textup{ for each integer } k, \mathbf x(k)\in \{0,2\}\big\}
$$
is used in the proof. We leave the details to a reader.
\end{enumerate}
\end{proof}
\begin{theorem}\label{cantorw}
The following hold for the sets of periodic points in $(X_F^+/_{\sim_+},(\sigma_F^+)^{\star})$ and in $(X_F/_{\sim},\sigma_F^{\star})$.
\begin{enumerate}
	\item\label{ejn} The set $\mathcal P((\sigma_F^+)^{\star})$ of periodic points in the quotient $(X_F^+/_{\sim_+},(\sigma_F^+)^{\star})$ is not dense in $X_F^+/_{\sim_+}$. 
	\item The set $\mathcal P(\sigma_F^{\star})$ of periodic points in the quotient $(X_F/_{\sim},\sigma_F^{\star})$ is not dense in $X_F/_{\sim}$. 
\end{enumerate}
	\end{theorem}
\begin{proof}
We prove each of the statements separately.
\begin{enumerate}
	\item Let $U=(0,1)\times (0,1)\times \prod_{k=3}^{\infty}X$. Then $U$ is open in $\prod_{k=1}^{\infty}X$ and $U\cap X_F^+\neq \emptyset$. However, note that $(U\cap \mathcal P(\sigma_F^+))=\emptyset$ (since for each $x\in (0,1)$, for each $\mathbf x=(x_1,x_2,x_3,\ldots)\in X_F^+$ such that $x_1=x$, and for each positive integer $n>1$, if $x_n\in (0,1)$, then there are positive integers $k$ and $\ell$ such that $x_n=x^{\frac{2^k}{3^{\ell}}}$, which is not equal to $x$). . It follows that the set $\mathcal P(\sigma_F^+)$ of periodic points in $(X_F^+,\sigma_F^+)$ is not dense in $X_F^+$. Therefore, by Theorem \ref{kvocki1}, the set $\mathcal P(\sigma_F^+)$ of periodic points in $(X_F^+,\sigma_F^+)$ is not dense in $X_F^+$.
	\item Suppose that the set $\mathcal P(\sigma_F^{\star})$ of periodic points in the quotient $(X_F/_{\sim},\sigma_F^{\star})$ is dense in $X_F/_{\sim}$. Therefore, by Theorem \ref{kvocki1}, the set $\mathcal P(\sigma_F)$ of periodic points in $(X_F,\sigma_F)$ is dense in $X_F$. It follows from  Theorem \ref{Mah}, the set $\mathcal P(\sigma^{-1})$ of periodic points in $(\varprojlim(X_F^+,\sigma_F^{+}),\sigma^{-1})$ is dense in $\varprojlim(X_F^+,\sigma_F^{+})$. By Theorem \ref{ingram1}, the set $\mathcal P(\sigma_F^+)$ of periodic points in $(X_F^+,\sigma_F^+)$ is dense in $X_F^+$, which contradicts with \ref{ejn}. 
\end{enumerate}
\end{proof}
\begin{theorem}\label{tttt}
	The dynamical systems $(X_F^+/_{\sim_+},(\sigma_F^+)^{\star})$ and $(X_F/_{\sim},\sigma_F^{\star})$ are both transitive.
	\end{theorem}
\begin{proof}
	The proof of this theorem is analogous to the proof of Theorem \ref{ttt}. We leave the details to a reader. 
\end{proof}

\begin{theorem}\label{robi}
The following hold for the Cantor fan $C$. 
\begin{enumerate}
\item There is a continuous mapping $f$ on the Cantor fan $C$, which is not a homeomorphism,  such that $(C,f)$ is mixing as well as chaotic in the sense of Robinson but not in the sense of Devaney.
	\item There is a homeomorphism $h$ on the Cantor fan $C$ such that $(C,h)$ is  mixing as well as chaotic in the sense of Robinson but not in the sense of Devaney. 
	\end{enumerate}
\end{theorem}
\begin{proof}
We prove each part of the theorem separately. 
\begin{enumerate}
	\item Let $C=X_F^+/_{\sim_+}$ and let $f=(\sigma_F^+)^{\star}$. Note that $f$ is a continuous function which is not a homeomorphism. By Theorem \ref{cantor}, $(C,f)$ has sensitive dependence on initial conditions. By Theorem \ref{tttt}, $(C,f)$ is transitive.  It follows from Theorem \ref{cantorw} that the set $\mathcal P(f)$ of periodic points in the quotient $(C,f)$ is not dense in $C$.  Therefore, $(C,f)$ is chaotic in the sense of Robinson but it is not chaotic in the sense of Devaney. 
	
	It follows from Theorem \ref{van} that $(X_F^+,\sigma_F^+)$ is mixing since $\Delta_X\subseteq F$. It follows from Theorem \ref{tupki} that $(C,f)$ is also mixing.
	\item Let $C=X_F/_{\sim}$ and let $h=\sigma_F^{\star}$. Note that $h$ is a homeomorphism. The rest of the proof is analogous to the proof above. We leave the details to a reader.
	 \end{enumerate}
 \end{proof}

   \subsection{Mixing and Knudsen's but not Devaney's chaos on the Cantor fan}\label{31} 

 Here, we study functions $f$ on the Cantor fan $C$ such that $(C,f)$ is mixing as we well as chaotic in the sense of Knudsen but not in the sense of Devaney.
\begin{definition}\label{nall1}
In this subsection, we use $X$ to denote $X=[0,1]\cup [2,3]\cup[4,5]$ and we let $f_1,f_2,f_3:X\rightarrow X$ to be the homeomorphisms from $X$ to $X$ that are defined by  
$$
f_1(x)=x,
$$
$$
f_2(x)=\begin{cases}
				x+2\text{;} & x\in [0,1]\\
				(x-2)^2\text{;} & x\in [2,3]\\
				x\text{;} & x\in [4,5]
			\end{cases} 
$$
$$
f_3(x)=\begin{cases}
				x\text{;} & x\in [0,1]\\
				x+2\text{;} & x\in [2,3]\\
				(x-4)^{\frac{1}{2}}+2\text{;} & x\in [4,5]			
				\end{cases}  
$$
for each $x\in X$. Then we use $F$ to denote the relation 
$$
F=\Gamma(f_1)\cup \Gamma(f_2)\cup \Gamma(f_3);
$$
 see Figure \ref{fig4}.
\begin{figure}[h!]
	\centering
		\includegraphics[width=17em]{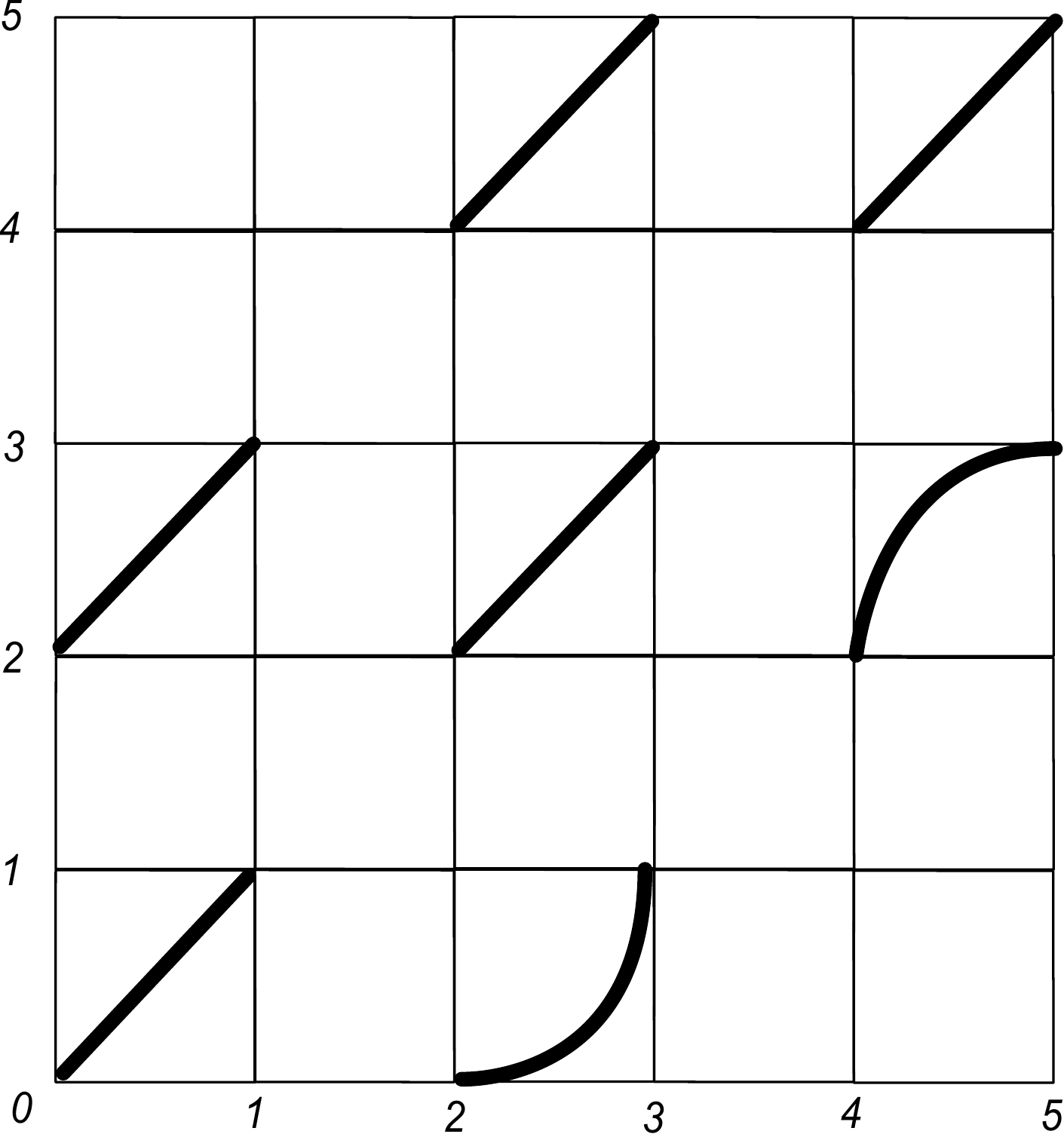}
	\caption{The  relation $F$ from Definition \ref{nall1}}
	\label{fig4}
\end{figure}  
	\end{definition}
	\begin{definition}
	We define two equivalence relations. 
	\begin{enumerate}
		\item For all $\mathbf x,\mathbf y\in X_F^+$, we define the relation $\sim_+$ as follows:
$$
\mathbf x\sim_+ \mathbf y ~~~   \Longleftrightarrow  ~~~  \mathbf x = \mathbf y \textup{ or  for each positive integer } k,  {\{\mathbf x(k),\mathbf y(k)\}\subseteq} \{0,2,4\}. 
$$
\item For all $\mathbf x,\mathbf y\in X_F$, we define the relation $\sim$ as follows:
$$
\mathbf x\sim \mathbf y ~~~   \Longleftrightarrow  ~~~  \mathbf x = \mathbf y \textup{ or  for each integer } k,  {\{\mathbf x(k),\mathbf y(k)\}\subseteq} \{0,2,4\}. 
$$
	\end{enumerate}
\end{definition}
\begin{observation}
	Essentially the same proof as the proof of \cite[Example 4.14]{BE} shows that the quotient spaces $X_F^+{/}_{\sim_+}$ and $X_F{/}_{\sim}$ are both  Cantor fans. Again, note that $(\sigma_F^+)^{\star}$ is not a homeomorphism on $X_F^+{/}_{\sim_+}$ while $\sigma_F^{\star}$ is a homeomorphism on $X_F{/}_{\sim}$.
	\end{observation}
	\begin{theorem}\label{skola0s}
		Let $A=\big\{\mathbf x\in X_F^+ \ | \ \textup{ for each positive integer } k, \mathbf x(k)\in \{0,2,4\}\big\}$.
		 Then
		 \begin{enumerate}
		 	\item $\sigma_F^+(A)\subseteq A$ and $\sigma_F^+(X_F^+\setminus A)\subseteq X_F^+\setminus A$, and
		 	\item $(X_F^+,\sigma_F^+)$ has sensitive dependence on initial conditions with respect to $A$.
		 \end{enumerate}  
	\end{theorem}
	\begin{proof}
		The proof is analogous to the proof of Theorem \ref{skola0}. We leave the details to a reader.
	\end{proof}	
	\begin{corollary}\label{corolars}
	Let $B=\big\{\mathbf x\in X_F^+ \ | \ \textup{ for each integer } k, \mathbf x(k)\in \{0,2,4\}\big\}$.
		 Then
		 \begin{enumerate}
		 	\item $\sigma_F(B)\subseteq B$ and $\sigma_F(X_F\setminus B)\subseteq X_F\setminus B$, and
		 	\item $(X_F,\sigma_F)$ has sensitive dependence on initial conditions with respect to $B$.
		 \end{enumerate}
	\end{corollary}
	\begin{proof}
	The proof is analogous to the proof of Corollary \ref{corolar}. We leave the details to a reader.
 		   	  \end{proof}
	
	\begin{theorem}\label{cantorss}
		The dynamical systems $(X_F^+/_{\sim_+},(\sigma_F^+)^{\star})$ and $(X_F/_{\sim},\sigma_F^{\star})$ both have sensitive dependence on initial conditions. 
	\end{theorem}
\begin{proof}
The proof is analogous to the proof of Theorem \ref{cantor}. We leave the details to a reader.
\end{proof}
\begin{theorem}\label{cantorwww}
The following hold for the sets of periodic points in $(X_F^+/_{\sim_+},(\sigma_F^+)^{\star})$ and $(X_F/_{\sim},\sigma_F^{\star})$.
\begin{enumerate}
	\item\label{ejn} The set $\mathcal P((\sigma_F^+)^{\star})$ of periodic points in the quotient $(X_F^+/_{\sim_+},(\sigma_F^+)^{\star})$ is dense in $X_F^+/_{\sim_+}$. 
	\item The set $\mathcal P(\sigma_F^{\star})$ of periodic points in the quotient $(X_F/_{\sim},\sigma_F^{\star})$ is dense in $X_F/_{\sim}$. 
\end{enumerate}
	\end{theorem}
\begin{proof}
The proof is analogous to the proof of Theorem \ref{ppp}. We leave the details to a reader.
\end{proof}
\begin{theorem}\label{transtrans}
The following hold for the dynamical systems $(X_F^+/_{\sim_+},(\sigma_F^+)^{\star})$ and $(X_F/_{\sim},\sigma_F^{\star})$.
\begin{enumerate}
	\item\label{ejn} The dynamical system $(X_F^+/_{\sim_+},(\sigma_F^+)^{\star})$ is not transitive. 
	\item The dynamical system $(X_F/_{\sim},\sigma_F^{\star})$ is not transitive.
	\end{enumerate}
	\end{theorem}
\begin{proof}
We prove each of the statements separately.
\begin{enumerate}
	\item To prove that $(X_F^+/_{\sim_+},(\sigma_F^+)^{\star})$ is not transitive, we show first that $(X_F^+,\sigma_F^+)$ is not transitive. Let $x\in X$. We consider the following cases. 
	\begin{enumerate}
	\item $x\in \{0,2,4\}$. Then $\mathcal U^{\oplus}_F(x)=\{0,2,4\}$.
	\item $x\in \{1,3,5\}$. Then $\mathcal U^{\oplus}_F(x)=\{1,3,5\}$.
	\item $x\not\in \{0,1,2,3,4,5\}$. Then 
	$$
	\mathcal U_F^{\oplus}(x)=\{x^{2^k}+2\ell \ | \ k \textup{ is an integer and } \ell\in \{0,1,2\}\}.
	$$
	\end{enumerate}
	In each case, $\mathcal U_F^{\oplus}(x)$ is not dense in $X$. For each $x\in X$, let $V_x$ be a non-empty open set in $X$ such that $V_x\cap \mathcal U_F^{\oplus}(x)=\emptyset$ and let $U_x=V_x\times \prod_{k=2}^{\infty}X$. 	It follows that for each $x\in X$ and for each point $\mathbf x=(x_1,x_2,x_3,\ldots)\in X_F^+$ such that $x_1=x$, 
	$$
	\{\mathbf x,\sigma_F^+(\mathbf x), (\sigma_F^+)^2(\mathbf x),(\sigma_F^+)^3(\mathbf x),\ldots\}\cap U_x=\emptyset.
	$$
	Therefore, for any $\mathbf x\in X_F^+$, the orbit $\{\mathbf x,\sigma_F^+(\mathbf x), (\sigma_F^+)^2(\mathbf x),(\sigma_F^+)^3(\mathbf x),\ldots\}$ of the point $\mathbf x$ is not dense in $X_F^+$. It follows from Theorem \ref{andrej} that $(X_F^+,\sigma_F^+)$ is not transitive. Therefore, by Proposition \ref{kvocienti}, the dynamical system $(X_F^+/_{\sim_+},(\sigma_F^+)^{\star})$ is not transitive. 
	\item Since $p_1(F)=p_2(F)=X$ and since $(X_F^+,\sigma_F^+)$ is not transitive, it follows from Theorem \ref{tazadnji} that the dynamical system  $(X_F,\sigma_F)$ is not transitive. Therefore, it follows from Theorem \ref{kvocienti} that the dynamical system $(X_F/_{\sim},\sigma_F^{\star})$ is not transitive. 
	\end{enumerate}
\end{proof}
\begin{theorem}\label{robicar}
The following hold for the Cantor fan $C$. 
\begin{enumerate}
\item There is a continuous mapping $f$ on the Cantor fan $C$, which is not a homeomorphism,  such that $(C,f)$ is mixing as well as chaotic in the sense of Knudsen but not in the sense of Devaney.
	\item There is a homeomorphism $h$ on the Cantor fan $C$ such that $(C,h)$ is mixing as well as chaotic in the sense of Knudsen but not in the sense of Devaney. 
	\end{enumerate}
\end{theorem}
\begin{proof}
We prove each part of the theorem separately. 
\begin{enumerate}
	\item Let $C=X_F^+/_{\sim_+}$ and let $f=(\sigma_F^+)^{\star}$. Note that $f$ is a continuous function which is not a homeomorphism. By Theorem \ref{cantorss}, $(C,f)$ has sensitive dependence on initial conditions. By Theorem \ref{transtrans}, $(C,f)$ is not transitive.  It follows from Theorem \ref{cantorwww} that the set $\mathcal P(f)$ of periodic points in  $(C,f)$ is dense in $C$.  Therefore, $(C,f)$ is chaotic in the sense of Knudsen but it is not chaotic in the sense of Devaney.
	
	It follows from Theorem \ref{van} that $(X_F^+,\sigma_F^+)$ is mixing since $\Delta_X\subseteq F$. It follows from Theorem \ref{tupki} that $(C,f)$ is also mixing.
	\item Let $C=X_F/_{\sim}$ and let $h=\sigma_F^{\star}$. Note that $h$ is a homeomorphism. The rest of the proof is analogous to the proof above. We leave the details to a reader.
	  \end{enumerate}
 \end{proof}

\section{Uncountable family of (non-)smooth fans that admit mixing homeomorphisms}\label{s5}
\begin{figure}[h!]
	\centering
		\includegraphics[width=30em]{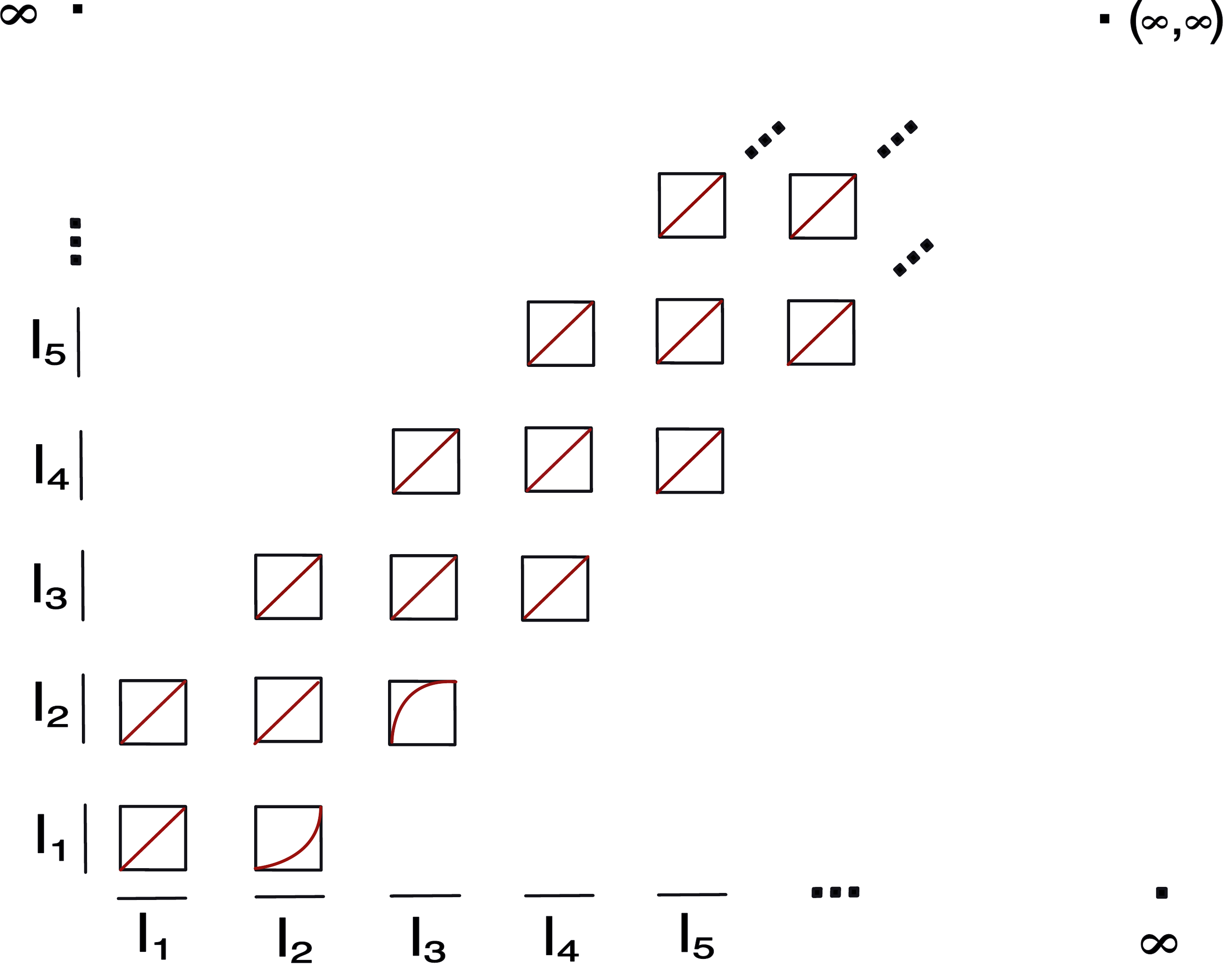}
	\caption{The relation $H$ on $X$}
	\label{uncun}
\end{figure}

In this section, an uncountable family $\mathcal G$ of pairwise non-homeomorphic smooth fans that admit mixing homeomorphisms is constructed. Our construction of the family $\mathcal G$ follows the idea from \cite{banic1}, where an uncountable family $\mathcal F$ of pairwise non-homeomorphic smooth fans that admit transitive homeomorphisms is constructed: every step of the construction of family $\mathcal F$ from \cite{banic1} is essentially copied here to construct the family $\mathcal G$. The only difference is a small modification of the relation $H$ on $X$ that is used in \cite{banic1} to obtain the family $\mathcal F$: in $H$, the graph in $(I_1\times I_1)\cup (I_2\times I_2)$ is replaced with the graph in $(I_2\times I_1)\cup (I_3\times I_2)$ and the graph in $(I_2\times I_1)\cup (I_3\times I_2)$ is replaced with the graph in $(I_1\times I_1)\cup (I_2\times I_2)$; see \cite[Figure 5]{banic1} and Figure \ref{uncun}. Therefore, in this section, we omit the details and simply state our first theorem.  
 
\begin{theorem}\label{mine}
There is a family $\mathcal G$ of uncountable many pairwise non-homeo\-morphic smooth fans that admit mixing homeomorphisms.
\end{theorem}

	In \cite{banic9}, a family of uncountably many pairwise non-homeomor\-phic non-smooth fans that admit transitive homeomorphisms is constructed from the family $\mathcal F$ from \cite{banic1}. This is done in such a way that for each smooth fan $F\in \mathcal F$, a special equivalence relation $\sim$ on $F$ is defined in such a way that $F/_{\sim}$ is a non-smooth fan that admits a transitive homeomorphism. The same procedure as the one from \cite{banic9} for the family $\mathcal F$, works also for our family $\mathcal G$. It transformes every smooth fan $F\in \mathcal G$ to a non-smooth fan $F/_{\sim}$ that admits a mixing homeomorphism. The following theorem is, therefore, a good place to finish the paper. Since its proof is essentially the same as the proof of \cite[Theorem 3.17]{banic9}, we leave the details to a reader.
	
\begin{theorem}\label{mine1}
There is a family of uncountable many pairwise non-homeo\-morphic non-smooth fans that admit mixing homeomorphisms.
\end{theorem}

\section{Acknowledgement}
This work is supported in part by the Slovenian Research Agency (research projects J1-4632, BI-HR/23-24-011, BI-US/22-24-086 and BI-US/22-24-094, and research program P1-0285). 
	

\noindent I. Bani\v c\\
              (1) Faculty of Natural Sciences and Mathematics, University of Maribor, Koro\v{s}ka 160, SI-2000 Maribor,
   Slovenia; \\(2) Institute of Mathematics, Physics and Mechanics, Jadranska 19, SI-1000 Ljubljana, 
   Slovenia; \\(3) Andrej Maru\v si\v c Institute, University of Primorska, Muzejski trg 2, SI-6000 Koper,
   Slovenia\\
             {iztok.banic@um.si}           
     
				\-
				
		\noindent G.  Erceg\\
             Faculty of Science, University of Split, Rudera Bo\v skovi\' ca 33, Split,  Croatia\\
{{goran.erceg@pmfst.hr}       }    

                 	\-
					
  \noindent J.  Kennedy\\
             Department of Mathematics,  Lamar University, 200 Lucas Building, P.O. Box 10047, Beaumont, Texas 77710 USA\\
{{kennedy9905@gmail.com}       }    

	\-
				
		\noindent C.  Mouron\\
             Rhodes College,  2000 North Parkway, Memphis, Tennessee 38112  USA\\\
{{mouronc@rhodes.edu}       }    

                 	\-
				
		\noindent V.  Nall\\
             Department of Mathematics,  University of Richmond, Richmond, Virginia 23173 USA\\
{{vnall@richmond.edu}       }   




\begin{thebibliography}{9}
\bibitem{akin} E.~Akin, General Topology of Dynamical Systems, Volume 1, Graduate Studies in Mathematics Series, American Mathematical Society, Providence RI, 1993.
\bibitem{aoki} N.~Aoki, Topological dynamics, in: K. Morita and J. Nagata, eds., Topics in General Topology (Elsevier, Amsterdam, 1989) 625--740.
\bibitem{EK} I.~Bani\v c, G.~Erceg,  J.~Kennedy, An embedding of the Cantor fan into the Lelek fan, https://web.math.pmf.unizg.hr/~duje/radhazumz/preprints/banic-erceg-kennedy-preprint.pdf
\bibitem{BE} I.~Bani\v c, G.~Erceg,   J.~Kennedy, C.~Mouron, V.~Nall, Transitive mappings on the Cantor fan,  
https://doi.org/10.48550/arXiv.2304.03350.
\bibitem{judyk} I.~Bani\v c, G.~Erceg,   J.~Kennedy,  V.~Nall, Quotients of dynamical systems and chaos on the Cantor fan, preprint.  
https://doi.org/10.48550/arXiv.2304.03350.
\bibitem{banic1} I.~Bani\v c, G.~Erceg,  J.~Kennedy, C.~Mouron, V.~Nall, An uncountable family of smooth fans that admit transitive homeomorphisms,  arXiv:2309.04003.
\bibitem{banic9} I.~Bani\v c,  J.~Kennedy, C.~Mouron, V.~Nall, An uncountable family of non-smooth fans that admit transitive homeomorphisms,  
https://doi.org/10.48550/arXiv.2310.08711.
\bibitem{banic2} I.~Bani\v c, G.~Erceg,  J.~Kennedy, A transitive homeomorphism on the Lelek fan,  {to appear in J. Difference Equ. Appl.  (2023) https://doi.org/10.1080/10236198.2023.2208242.}
\bibitem{banks} J. ~Banks, J. ~Brooks, G. ~Cairns, G. ~Davis and P. ~Stacey, On Devaney's Definition of Chaos, The American Mathematical Monthly 99 (1992) 332--334.
\bibitem{munkres} J. R. Munkres, Topology: a first course, Prentice-Hall, Inc., Englewood Cliffs, N.J., 1975
\bibitem{oversteegen} W.~D.~Bula and L.~Overseegen, A Characterization of smooth Cantor Bouquets,  Proc. Amer.Math.Soc. 108 (1990) 529--534.
\bibitem{Jcharatonik}J.~J. ~Charatonik,  On fans,  Dissertationes Math.  54 (1967).
\bibitem{charatonik} W.~J.~Charatonik, The Lelek fan is unique, Houston J. Math. 15 (1989) 27--34.
\bibitem{charatonik3} J. ~J.~ Charatonik, On fans. Dissertationes Math. (Rozprawy Mat.) 7133 (1967), 37 pp.
\bibitem{devaney} R.~L.~Devaney, A first course in chaotic dynamical systems: theory and experiments. Massachusetts: Perseus Books, 1992.
\bibitem{eberhart} C. ~Eberhart, A note on smooth fans, Colloq.  Math.  20 (1969) 89--90.
\bibitem{engelking1} R.~ Engelking, General topology, Heldermann, Berlin, 1989.
\bibitem{KS} S.~Kolyada, L.~Snoha,  Topological transitivity, \textit{Scholarpedia}  4 (2):5802 (2009).
\bibitem{knudsen} C.~Knudsen, Chaos Without Nonperiodicity, The American Mathematical Monthly 101(1994) 563--565.
\bibitem{koch} R.~J.~Koch, Arcs in partially ordered spaces, Pacific J.  Math.  20 (1959) 723--728
\bibitem{lelek} A.~Lelek, On plane dendroids and their end-points in the classical sense, Fund. Math. 49 (1960/1961) 301--319.
\bibitem{nadler} S.~B.~Nadler, Continuum theory. An introduction, Marcel Dekker, Inc., New York, 1992.
\bibitem{robinson} C.~Robinson,  Dynamical systems: stability, symbolic dynamics, and chaos. 2nd ed. Boca Raton, FL: CRC Press Inc., 1999.
\bibitem{willard} S.~ Willard, General topology, Dover Publications, New York, 1998.
\bibitem{wu} Xinxing Wu, Xiong Wang, Guanrong Chen, $\mathcal F$-mixing property and $(\mathcal F_1, \mathcal F_2)$-everywhere chaos of inverse limit dynamical systems, Nonlinear Analysis: Theory, Methods and Applications 104 (2014) 147--155.
\end{thebibliography}
\end{document}